\documentclass[%
 amssymb, amsmath,%
 aip,jmp
,preprint%
]{revtex4-1}

\usepackage{docs}%
\usepackage{bm}%
\usepackage{float}
\usepackage[colorlinks=true,linkcolor=blue]{hyperref}%
\usepackage{caption}
\DeclareCaptionType{algorithm}
\usepackage{algorithmic}
\expandafter\ifx\csname package@font\endcsname\relax\else
 \expandafter\expandafter
 \expandafter\usepackage
 \expandafter\expandafter
 \expandafter{\csname package@font\endcsname}%
\fi
\def\D{\displaystyle}

\DeclareMathOperator{\ord}{ord}
\DeclareMathOperator{\lt}{lt}
\DeclareMathOperator{\lc}{lc}
\DeclareMathOperator{\lcm}{lcm}
\DeclareMathOperator{\lex}{lex}

\begin{document}

\allowdisplaybreaks
\def\D{\displaystyle}
\newtheorem{theorem}{Theorem}[section]
\newtheorem{definition}[theorem]{Definition}
\newtheorem{lemma}[theorem]{Lemma}
\newtheorem{example}[theorem]{Example}
\newtheorem{Diffusionexample}[theorem]{Example (Diffusion equation in $1$-space)}
\newtheorem{Maxwellexample}[theorem]{Example (Maxwell equations for vanishing free current density and \newline free charge density)}
\newtheorem{proposition}[theorem]{Proposition}
\newtheorem{remark}[theorem]{Remark}
\newtheorem{note}[theorem]{Note}
\newtheorem{deflemma}[theorem]{Lemma and Definition}

\title[Computation of the Strength]{Computation of the Strength of PDEs of Mathematical Physics and their Difference Approximations}%

\author{Christian D\"onch}
\email{cdoench@risc.jku.at}
\affiliation{Research Institute for Symbolic Computation, Johannes Kepler University Linz\\
 A-4040 Linz, Austria}
\author{Alexander Levin}%
\email{levin@cua.edu}
\affiliation{The Catholic University of America\\
Washington, D. C. 20064, USA}%

\date{February 2012}%
\begin{abstract}
We develop a method for evaluation of A. Einstein's strength of systems of partial differential and difference equations based on the computation of Hilbert-type dimension polynomials of the associated differential and difference field extensions. Also we present algorithms for such computations, which  are based on the Gr\"obner basis method adjusted for the modules over rings of differential, difference and inversive difference operators. The developed technique is applied to some fundamental systems of PDEs of mathematical physics such as the diffusion equation, Maxwell equations and equations for an electromagnetic field given by its potential. In each of these cases we determine the strength of the original system of PDEs and the strength of the corresponding systems of partial difference equations obtained by forward and symmetric difference schemes. In particular, we obtain a method for comparing two difference schemes from the point of view of their strength.
\end{abstract}

\maketitle


\section{Introduction}

The concept of the strength of a system of partial differential equations (PDEs) was introduced by A. Einstein  as a measure for the size of the solution space of such a system. In \cite{Einstein} A. Einstein defined the strength of a system of partial differential equations governing a physical field as follows: ''\ldots the system of equations is to be chosen so that the field quantities are determined as strongly as possible. In order to apply this principle, we propose a method which gives a measure of strength of an equation system. We expand the field variables, in the neighborhood of a point $\mathcal{P}$, into a Taylor series (which presupposes the analytic character of the field); the coefficients of these series, which are the derivatives of the field variables at $\mathcal{P}$, fall into sets according to the degree of differentiation. In every such degree there appear, for the first time, a set of coefficients which would be free for arbitrary choice if it were not that the field must satisfy a system of differential
equations.  Through this system of differential equations (and its derivatives with respect to the coordinates) the number of coefficients is restricted, so that in each degree a smaller number of coefficients is left free for arbitrary choice. The set of numbers of 'free' coefficients for all degrees of differentiation is then a measure of the 'weakness' of the system of equations, and through this, also of its 'strength'.''

Calculating by hand A. Einstein found out, that, for example, the potential and field formulations of Maxwell equations have different strengths for the dimension four. However, he did not obtain the exact expression of the above-mentioned number of free coefficients as a function of the degree of differentiation. Even though there were a number of works on the strength of a system of differential equations (in particular, on its relation to Cartan characters), see, for example, \cite{Mariwalla,Matthews1,Matthews2,Schutz,Seiler1,Seiler2,Seiler3}, and  \cite{Sue}, there was no method of evaluating such a function until 1980 when A. Mikhalev and E. Pankratev \cite{MP}  showed that the strength of a system of algebraic partial differential equations (that is, a system of the form $f_{i} = 0$, $i\in I$, where $f_{i}$ are multivariate polynomials in unknown functions and their partial derivatives) is expressed by Kolchin's differential dimension polynomial associated with the differential field extension defined by the system.  This observation allowed A. Mikhalev and E. Pankratev  to develop two methods of determining the strength of a system of algebraic PDEs via computing the differential dimension polynomial of the corresponding differential field extension.  The first method is based on construction of a characteristic set of the ideal of differential polynomials defined by the system and then computing the differential dimension polynomial using the leading terms of the elements of the characteristic set (the idea of this approach comes from the original proof of Kolchin's theorem, see \cite[Chapter II, Theorem 6]{Kolchin}). The second approach is based on the works by J. Johnson \cite{Johnson1,Johnson2}, who showed that the differential dimension polynomial of a differential field extension can be computed as a Hilbert polynomial of the associated module of K\"ahler differentials. Using free resolutions for such a module,  A. Mikhalev and E. Pankratev \cite{MP} evaluated the strength of several well-known systems of PDEs including the wave equation, both forms of Maxwell equations, Dirac equations (with zero mass), Lame equations, and some other systems of PDEs of mathematical physics. Note that A. Einstein, K. Mariwalla, M. Sue and some other authors who investigated the concept of strength in 1970s characterized the strength of a system by the ''coefficient of freedom'', an integer, that is fully determined by the leading coefficient of the differential dimension polynomial. The fact that such a polynomial provides a far more precise description of the strength than its leading term was justified by the result of W. Sit \cite{Sit} who proved that the set of differential dimension polynomials is well-ordered with respect to the natural order ($f(t) < g(t)$ if and only if $f(r) < g(r)$ for all sufficiently large integers $r$); this result  allows one to distinguish two systems of PDEs with the same ''coefficient of freedom'' by their strength.

Since 1980s the technique of dimension polynomials has been extended to the analysis of systems of algebraic difference and difference-differential equations. In a series of works whose results are summarized in \cite{Levin7} the second author proved the existence and developed some methods of computation of dimension polynomials of difference field extensions and systems of algebraic difference equations. These polynomials determine A. Einstein's strength of a system of algebraic partial difference equations (we give the details in Section 3 of this work) and, in particular,  allow one to evaluate the quality of  difference schemes for PDEs from the point of view of their strength.

The next step in the analysis of systems of  PDEs and systems of partial difference equations is to consider their degrees of freedom with respect to different groups of basic operators (differentiations or translations). Theorems on multivariate dimension polynomials proved in \cite{Levin4,Levin5,Levin6} (see also \cite[Chapters 3, 4, 7]{Levin7}) allow to characterize the strength of a system of partial differential, difference or difference-differential equations in the case when the ''weights'' of basic operators  of different groups are different. Methods of computations of multivariate dimension polynomials for systems of differential, difference and difference-differential equations were developed in \cite{Levin4,Levin5,Levin6,Levin7,Donch,DW,ZW1}, and \cite{ZW2} with the use of generalizations of the Gr\"obner basis technique. In particular, the first author has implemented in Maple two algorithms of computation of bivariate difference-differential dimension polynomials via relative Gr\"obner bases introduced in \cite{ZW1}.

In this paper we present the theory and technique of differential, difference, and difference-differential dimension polynomials together with the applications of this technique to the analysis of fundamental systems of PDEs of mathematical physics and corresponding systems of partial difference equations. In particular, we develop a method that allows one to compute the strength of such systems in the sense of A. Einstein and compare different difference schemes for a given system of PDEs by their strength. We illustrate this method with the computation of the strength of the diffusion equation, Maxwell equations and equations for an electromagnetic field given by its potential, as well as with the computation of the strength of systems of difference equations obtained from these PDEs via different difference schemes.

\section{Preliminaries}

In this section we present some basic concepts and results that are used throughout the paper. In what follows, ${\mathbb N}, {\mathbb Z}$,
${\mathbb Q}$, and ${\mathbb R}$ denote the sets of all non-negative integers, integers, rational numbers, and real numbers, respectively. The number of elements of a set $A$ is denoted by $|A|$. As usual, ${\mathbb Q}[t]$ denotes the ring of polynomials in one variable $t$ with rational coefficients. By a ring we always mean an associative ring with unit element. Every ring homomorphism is unitary (maps unit element onto unit element), every subring of a ring contains the unit element of the ring. Unless otherwise indicated, by a module over a ring $R$ we always mean a unitary left $R$-module.

\medskip

{\bf 2.1.\, Differential and difference rings and fields}\,

\medskip

A {\em differential ring} (respectively, a {\em difference ring}) is a commutative ring $R$ together with a finite set $\Delta = \{\delta_{1}, \dots, \delta_{m}\}$ of mutually commuting mappings of $R$ into itself such that each $\delta_{i}$ is a derivation of $R$ (respectively, $\delta_{i}$ are  injective endomorphisms of $R$ also called {\em translations}). The set $\Delta$ is said to be the {\em basic set\/} of the differential (or difference) ring $R$, which is also called a {\em $\Delta$-ring}. If a $\Delta$-ring is a field, it is called a $\Delta$-field (this is a differential field if $\Delta$ is a set of mutually commuting derivations and a difference field if the elements of $\Delta$ are endomorphisms). If $\delta_{1}, \dots, \delta_{m}$ are automorphisms of a difference ring $R$, we say that $R$ is an {\em inversive difference ring\/} with the basic set $\Delta$.  In this case we denote the set $\{\delta_{1}, \dots, \delta_{m}, \delta_{1}^{-1}, \dots, \delta_{m}^{-1}\}$ by $\Delta^{\ast}$ and call $R$ a $\Delta^{\ast}$-ring (if $R$ is a field, it is called an inversive difference field or a $\Delta^{\ast}$-field).

Let $R$ be a $\Delta$-ring and $R_{0}$ a subring (ideal) of $R$ such that $\delta(R_{0})\subseteq R_{0}$ for any $\delta\in\Delta$. Then $R_{0}$ is
called a $\Delta$-subring (respectively, a $\Delta$-ideal)  of $R$. If $R_{0}$ is a $\Delta$-subring of $R$, we also say that $R$ is a $\Delta$-ring extension of $R_{0}$.

If the elements of $\Delta$ act on $R$ as mutually commuting derivations, we say that $R_{0}$ is a differential subring (differential ideal) of $R$; if the elements of $\Delta$ are mutually commuting injective endomorphisms, we say that $R_{0}$ is a difference subring (difference ideal) of $R$. Anyway, the prefix $\Delta$-, depending on the context,  means either ''differential'' or ''difference'', while the prefix $\Delta^{\ast}$- means ''inversive difference''.   If $R$ is a $\Delta^{\ast}$-ring (this assumption implies that $\Delta$ is a set of mutually commuting automorphisms of $R$), then a subring (ideal) $R_{0}$ of $R$ is called a $\Delta^{\ast}$-subring (respectively, $\Delta^{\ast}$-ideal) of $R$ if $\alpha(R_{0})\subseteq R_{0}$ for any $\alpha\in \Delta^{\ast}$. ($\Delta^{\ast}$-ideals are also called reflexive difference ideals of $R$; this term, as well as the term $\Delta^{\ast}$-ideal, is also used for ideals $I$ of a difference $\Delta$-ring $R$ such that for any $\delta\in \Delta,\, a\in R$, the inclusion $\delta(a)\in I$ implies $a\in I$). If $R$ is a $\Delta$- (or $\Delta^{\ast}$-) field and $R_{0}$ a subfield of $R$ which is also a $\Delta$- (respectively, $\Delta^{\ast}$-) subring of $R$, then  $R_{0}$ is said to be a $\Delta$- (respectively, $\Delta^{\ast}$-) subfield of $R$ while $R$ is called a $\Delta$- (respectively, $\Delta^{\ast}$-) field extension (or overfield) of $R_{0}$. In this case we also say that we have a $\Delta$-(or $\Delta^{\ast}$-) field extension $R/R_{0}$.

\medskip

If $R$ is a $\Delta$-ring with a basic set $\Delta = \{\delta_{1}, \dots, \delta_{m}\}$, then $\Theta_{\Delta}$ (or $\Theta$ if the set $\Delta$
is fixed) will denote the free commutative monoid generated by $\delta_{1}, \dots, \delta_{m}$. Elements of $\Theta$ will be written in the multiplicative form $\delta_{1}^{k_{1}}\dots \delta_{m}^{k_{m}}$ ($k_{1},\dots, k_{m}\in {\mathbb N}$) and considered as the corresponding mappings of $R$ into itself. If $R$ is an inversive difference ($\Delta^{\ast}$-) ring, then $\Gamma_{\Delta}$ (or $\Gamma$ if the set $\Delta$ is fixed) will denote the free commutative group generated by the set $\Delta$. It is clear that elements of the group $\Gamma$ (written in
the multiplicative form $\delta_{1}^{i_{1}}\dots \delta_{m}^{i_{m}}$ where $i_{1},\dots, i_{m}\in {\mathbb Z}$) act on $R$ as automorphisms
and $\Theta$ is a subsemigroup of $\Gamma$. 

Let $R$ be a $\Delta$-ring and $S\subseteq R$. Then the intersection of all $\Delta$-ideals of $R$ containing $S$ is denoted by $[S]$. Clearly, $[S]$ is the smallest $\Delta$-ideal of $R$ containing $S$; as an ideal, it is generated by the set  $\Theta S = \{\theta(a) | \theta \in \Theta, \,a\in
S\}$. If $J = [S]$, we say that the $\Delta$-ideal $J$ is generated by the set $S$ called a {\em set of $\Delta$-generators} of $J$. If $S$ is finite, $S=\{a_{1},\dots, a_{k}\}$, we write $J = [a_{1},\dots, a_{k}]$ and say that $J$ is a {\em finitely generated $\Delta$-ideal} of the $\Delta$-ring $R$. (In this case elements $a_{1},\dots, a_{k}$ are said to be {\em $\Delta$-generators} of $J$.) If $R$ is an inversive difference ($\Delta^{\ast}$-) ring and $S\subseteq R$, then the smallest $\Delta^{\ast}$-ideal of $R$ containing $S$ is denoted by $[S]^{\ast}$ (as an ideal, it is generated by the set $\Gamma S = \{\gamma(a) | \gamma \in \Gamma,\, a\in S\}$. If $S$ is finite, $S=\{a_{1},\dots, a_{k}\}$, we write $[a_{1},\dots, a_{k}]^{\ast}$ for $I = [S]^{\ast}$ and say that $I$ is a {\em finitely generated $\Delta^{\ast}$-ideal} of $R$; in this case the elements $a_{1},\dots, a_{k}$ are called  $\Delta^{\ast}$-generators of $I$.

\medskip

Let $R$ be a $\Delta$-ring, $R_{0}$ a $\Delta$-subring of $R$ and $B\subseteq R$. The intersection of all $\Delta$-subrings of $R$ containing $R_{0}$ and $B$ is called the {\em $\Delta$-subring of $R$ generated by the set $B$ over $R_{0}$}; it is denoted by $R_{0}\{B\}$. (As a ring, $R_{0}\{B\}$ coincides with the ring $R_{0}[\{\theta(b) | b\in B, \theta\in \Theta\}]$ obtained by adjoining the set $\{\theta(b) | b\in B, \theta \in \Theta\}$ to the ring $R_{0}$). The set $B$ is said to be the set of {\em $\Delta$-generators} of the $\Delta$-ring $R_{0}\{B\}$ over $R_{0}$. If this set is finite, $B = \{b_{1},\dots, b_{k}\}$, we say that $R' = R_{0}\{B\}$ is a finitely generated $\Delta$-ring extension (or $\Delta$-ring extension) of $R_{0}$ and write $R' = R_{0}\{b_{1},\dots, b_{k}\}$. If $R$ is a $\Delta$-field, $R_{0}$ a $\Delta$-subfield of $R$ and $B\subseteq R$, then the intersection
of all $\Delta$-subfields of $R$ containing $R_{0}$ and $B$ is denoted by $R_{0}\langle B\rangle$ (or $R_{0}\langle b_{1},\dots, b_{k}\rangle$ if $B=\{b_{1},\dots, b_{k}\}$ is a finite set). This is the smallest $\Delta$-subfield of $R$ containing $R_{0}$ and $B$; it coincides with the field $R_{0}(\{\theta(b) | b\in B, \theta\in\Theta\})$. The set $B$ is called a set of {\em $\Delta$-generators} of the $\Delta$-field $R_{0}\langle B\rangle$ over $R_{0}$. If $R$ is an inversive difference ($\Delta^{\ast}$-) ring, $R_{0}$ a $\Delta^{\ast}$-subring of $R$ and $B\subseteq R$. Then
the intersection of all $\Delta^{\ast}$-subrings of $R$ containing $R_{0}$ and $B$ is the smallest $\Delta^{\ast}$-subring of $R$ containing $R_{0}$ and $B$. This ring coincides with the ring $R_{0}[\{\gamma(b) | b\in B, \gamma \in \Gamma\}]$; it is denoted by $R_{0}\{B\}^{\ast}$. The set $B$ is said to be a {\em set of $\Delta^{\ast}$-generators} of $R_{0}\{B\}^{\ast}$ over $R_{0}$. If $B = \{b_{1},\dots, b_{k}\}$ is a finite set, we say that $S = R_{0}\{B\}^{\ast}$ is a finitely generated $\Delta^{\ast}$-) ring extension (or $\Delta^{\ast}$-overring) of $R_{0}$ and write $S = R_{0}\{b_{1},\dots, b_{k}\}^{\ast}$. If $R$ is a $\Delta^{\ast}$-field, $R_{0}$ a $\Delta^{\ast}$-subfield of $R$ and $B\subseteq R$, then the
intersection of all $\Delta^{\ast}$-subfields of $R$ containing $R_{0}$ and $B$ is denoted by $R_{0}\langle B\rangle^{\ast}$. This field
coincides with the field $R_{0}(\{\gamma(b) | b\in B, \gamma \in \Gamma\})$. The set $B$ is called a {\em set of $\Delta^{\ast}$-generators of the $\Delta^{\ast}$-field extension} $R_{0}\langle B\rangle^{\ast}$ {\em of} $R_{0}$.  If $B$ is finite, $B=\{b_{1},\dots, b_{k}\}$, we write $R_{0}\langle b_{1},\dots, b_{k}\rangle^{\ast}$ for $R_{0}\langle B\rangle^{\ast}$.

In what follows we often consider two or more $\Delta$- (or $\Delta^{\ast}$-) rings $R_{1},\dots, R_{p}$ with the same basic set $\Delta=\{\delta_{1},\dots, \delta_{m}\}$. Formally speaking, it means that for every $i=1,\dots, p$, there is some fixed mapping $\nu_{i}$ from the set $\Delta$ into the set of all derivations or injective endomorphisms of the ring $R_{i}$ such that any two mappings $\nu_{i}(\delta_{j})$ and $\nu_{i}(\delta_{k})$ of $R_{i}$ commute ($1\leq j, k\leq n$). We shall identify elements $\delta_{j}$ with their images $\nu_{i}(\delta_{j})$ and say that elements of the set $\Delta$ act as mutually commuting derivations or injective endomorphisms of the ring $R_{i}$ ($i=1,\dots, p$).

Let $R_{1}$ and $R_{2}$ be differential or difference or inversive difference rings with the same basic set $\Delta$.  A ring homomorphism
$\phi: R_{1}\rightarrow R_{2}$ is called a $\Delta$-homomorphism if $\phi(\delta(a)) = \delta(\phi(a))$ for any $\delta \in \Delta, a\in R_{1}$. (Clearly, if $\phi: R_{1} \rightarrow R_{2}$ is a $\Delta$-homomorphism of $\Delta^{\ast}$-rings, then $\phi(\delta(a)) = \delta(\phi(a))$ for any $\delta \in \Delta^{\ast}, a\in R_{1}$.)  If $R_{1}$ and $R_{2}$ are two $\Delta$-overrings of the same $\Delta$-ring $R_{0}$ and
$\phi: R_{1}\rightarrow R_{2}$ is a $\Delta$-homomorphism such that $\phi(a) = a$ for any $a\in R_{0}$, we say that $\phi$ is a $\Delta$-homomorphism over $R_{0}$ or that $\phi$ leaves the ring $R_{0}$ fixed. It is easy to see that the kernel of any $\Delta$-homomorphism of $\Delta$-rings $\phi: R\rightarrow R'$ is a $\Delta$-ideal of $R$ (moreover, in the case of difference rings, this kernel is a reflexive difference ideal of $R_{1}$). Conversely, let $g$ be a surjective homomorphism of a $\Delta$-ring $R$ onto a ring $S$ such that Ker$\,g$ is a $\Delta$- or $\Delta^{\ast}$- (if $R$ is a difference $\Delta$-ring) ideal of $R$. Then there is a unique structure of a $\Delta$-ring on $S$ such that $g$ is a $\Delta$-homomorphism. In particular, if $I$ is a $\Delta$- or $\Delta^{\ast}$- (if $R$ is a difference $\Delta$-ring) ideal of a $\Delta$-ring $R$, then the factor ring $R/I$ has a unique structure of a $\Delta$-ring such that the canonical surjection $R\rightarrow R/I$ is a $\Delta$-homomorphism.
In this case $R/I$ is said to be the $\Delta$-factor ring of $R$ by $I$.

\medskip

If a $\Delta$- (or $\Delta^{\ast}$-) ring $R$ is an integral domain, then its quotient field $Q(R)$ can be naturally considered as a $\Delta$- (respectively, $\Delta^{\ast}$-) overring of $R$. (If $\Delta$ consists of derivations, then they extend to $Q(R)$ via the quotient rule).  In this case $Q(R)$ is said to be the {\em quotient} $\Delta$- (respectively, $\Delta^{\ast}$-) {\em field} of $R$. Clearly, if a $\Delta$- (or $\Delta^{\ast}$-) field $K$ contains an integral domain $R$ as a $\Delta$- (respectively, $\Delta^{\ast}$-) subring, then $K$ contains the quotient $\Delta$- (respectively, $\Delta^{\ast}$-) field $Q(R)$.

\bigskip

{\bf 2.2.\, Differential, difference, and inversive difference polynomials. Algebraic differential and difference equations.}

\medskip

\, With the above notation, let $R$ be a $\Delta$- (or $\Delta^{\ast}$-) ring and let $U = \{u_{\lambda} | \lambda \in \Lambda\}$ be a family of elements in some $\Delta$- (respectively, $\Delta^{\ast}$-) ring extension of $R$. We say that the family $U$ is $\Delta$-{\em algebraically dependent} over $R$, if the family $\Theta U = \{\theta(u_{\lambda}) | \theta \in \Theta, \lambda \in \Lambda\}$ is algebraically dependent over $R$ (that is, there exist elements $v_{1},\dots, v_{k}\in \Theta U$ and a non-zero polynomial $f(X_{1},\dots, X_{k})$ with coefficients in $R$ such that $f(v_{1},\dots, v_{k}) = 0$). Otherwise, the family $U$ is said to be $\Delta$-{\em algebraically independent} over $R$ or a family of $\Delta$-{\em indeterminates} over $R$.   In the last case, the $\Delta$-ring $S=R\{(u_{\lambda})_{\lambda \in \Lambda}\}$ is called the {\em algebra of $\Delta$-polynomials} over $R$. If $\Delta$ consists of derivations (respectively, injective endomorphisms), then $S$ is also called a ring of differential (respectively, difference)  polynomials in the difference (or $\sigma$-) indeterminates $\{(u_{\lambda})_{\lambda \in \Lambda}\}$ over $R$. If $R$ is a $\Delta^{\ast}$-ring and the family $U$ considered above is $\Delta$-algebraically independent over $R$, then the ring $R\{(u_{\lambda})_{\lambda \in \Lambda}\}^{\ast}$ is called the algebra of $\Delta^{\ast}$-polynomials in the $\Delta^{\ast}$-indeterminates $u_{\lambda}$ over $R$.

If a family consisting of one element $u$ is $\Delta$-algebraically dependent over $R$, the element $u$ is said to be {\em $\Delta$-algebraic} over $R$. If the set  $\{\theta(u) | \theta\in\Theta\}$ is algebraically independent over $R$, we say that $u$ is $\Delta$-{\em transcendental} over the ring $R$.

Let $R$ be a $\Delta$-field, $L$ a $\Delta$-field extension of $R$, and $A\subseteq L$. We say that the set $A$ is {\em $\Delta$-algebraic over $R$} if every element $a\in A$ is $\Delta$-algebraic over $R$. If every element of $L$ is $\Delta$-algebraic over $R$, we say that $L$ is a {\em $\Delta$-algebraic field extension} of $R$.

 The following statement is proved in \cite[Chapter 1, Section 6]{Kolchin}, \cite[Chapter 2, Theorem I]{Cohn}, and \cite[Propositions 3.3.7, 3.4.4]{KLMP} for differential, difference and inversive difference rings.

\begin{proposition} .  Let $R$ be a $\Delta$- (respectively, $\Delta^{\ast}$-) ring and $I$ an arbitrary set. Then there exists an algebra of $\Delta$- (respectively, $\Delta^{\ast}$-) polynomials over $R$ in a family of $\Delta$- (respectively, $\Delta^{\ast}$-) indeterminates with indices from the set $I$. If $S$ and $S'$ are two such algebras, then there exists a $\Delta$-isomorphism $S\rightarrow S'$ that leaves the ring $R$ fixed. If $R$ is an integral domain, then any algebra of $\Delta$- (respectively, $\Delta^{\ast}$-) polynomials over $R$ is an integral domain.
\end{proposition}

The algebra of $\Delta$-polynomials over a $\Delta$-ring $R$ in a family of $\Delta$-indetermi\-nates with indices from a set $I$ is a polynomial $R$-algebra in the set of indeterminates $\Theta Y= \{y_{i,\theta}\}_{i\in I, \theta \in \Theta}$ with indices from the set $I\times \Theta$. This algebra, as it is shown in \cite[Chapter 1, Section 6]{Kolchin}, \cite[Chapter 2, Theorem I]{Cohn}, and \cite[Propositions 3.3.7]{KLMP} can be viewed as a $\Delta$-ring extension of $R$ where $\delta(y_{i,\theta}) = y_{i, \delta\theta}$ for any $\delta\in \Delta,\, y_{i,\theta}\in \Theta Y$. Setting $y_{i}=y_{i,1}$ we can write $y_{i,\theta}$ as $\theta y_{i}$. If $R$ is a $\Delta^{\ast}$-ring, then the algebra of $\Delta^{\ast}$-polynomials over $R$ in $\Delta$-indeterminates with indices from a set $I$ is a polynomial $R$-algebra $S$ in the set of indeterminates $\Gamma Y= \{y_{i,\gamma}\}_{i\in I, \gamma \in \Gamma}$ with indices from the set $I\times \Gamma$. As it is shown in \cite[Propositions 3.4.4]{KLMP}, $S$ can be treated as a $\Delta^{\ast}$-ring extension of $R$ where $\delta(y_{i,\gamma}) = y_{i, \delta\gamma}$ for any $\delta\in \Delta^{\ast},\, y_{i,\gamma}\in \Gamma Y$. In what follows we denote $y_{i,1}$ by $y_{i}$ and write $y_{i,\gamma}$ ($\gamma\in \Gamma$) as $\gamma y_{i}$.

\medskip

Let $R$ be a $\Delta$-ring, $R\{(y_{i})_{i\in I}\}$ an algebra of difference polynomials in a family of $\Delta$-indeterminates $\{(y_{i})_{i\in I}\}$, and $\{(\eta_{i})_{i\in I}\}$ a set of elements in some $\Delta$-ring extension of $R$. Since the set $\{\theta_{i}) | i\in I, \theta\in \Theta\}$ is algebraically independent over $R$, there exists a unique ring homomorphism $\phi_{\eta}: R[(\theta y_{i})_{i\in I, \theta \in \Theta}]\rightarrow
R[\theta(\eta_{i})_{i\in I, \theta \in\Theta}]$ that maps every $\theta y_{i}$ onto $\theta(\eta_{i})$ and leaves $R$ fixed. Clearly, $\phi_{\eta}$ is a surjective $\Delta$-homomorphism of $R\{(y_{i})_{i\in I}\}$ onto $R\{(\eta_{i})_{i\in I}\}$; it is called the {\em substitution} of $(\eta_{i})_{i\in I}$ for $(y_{i})_{i\in I}$. Similarly, if $R$ is a $\Delta^{\ast}$-ring, $R\{(y_{i})_{i\in I}\}^{\ast}$ an algebra of $\Delta^{\ast}$-polynomials over $R$ and $(\eta_{i})_{i\in I}$ a family of elements in a $\Delta^{\ast}$-ring extension of $R$, one can define a surjective $\Delta$-homomorphism $R\{(y_{i})_{i\in I}\}^{\ast}\rightarrow R\{(\eta_{i})_{i\in I}\}^{\ast}$ that maps every $y_{i}$ onto $\eta_{i}$ and leaves the ring $R$ fixed. This homomorphism is also called the substitution of $(\eta_{i})_{i\in I}$ for $(y_{i})_{i\in I}$. (It will be always clear whether we talk about substitutions for difference ($\Delta$-) or inversive difference ($\Delta^{\ast}$-) polynomials.) If $g$ is a $\Delta$- or $\Delta^{\ast}$- polynomial, then its image under a substitution of $(\eta_{i})_{i\in I}$ for $(y_{i})_{i\in I}$ is denoted by $g((\eta_{i})_{i\in I})$. The kernel of a   substitution is a $\Delta$- (or $\Delta^{\ast}$- if we deal with difference or inversive difference polynomials) ideal of the $\Delta$-ring $R\{(y_{i})_{i\in I}\}$ (respectively, of the $\Delta^{\ast}$-ring $R\{(y_{i})_{i\in I}\}^{\ast}$ if we consider substitution for inversive difference polynomials). This kernel is called the {\em defining} $\Delta$- (or $\Delta^{\ast}$-) {\em ideal} of the family
$(\eta_{i})_{i\in I}$ over $R$.

\medskip

If $R$ is a $\Delta$- (or $\Delta^{\ast}$-) field and $(\eta_{i})_{i\in I}$ is a family of elements in some $\Delta$- (respectively, $\Delta^{\ast}$-) overfield $S$, then $R\{(\eta_{i})_{i\in I}\}$ (respectively, $R\{(\eta_{i})_{i\in I}\}^{\ast}$) is an integral domain (it is
contained in the field $S$). It follows that the defining $\Delta$- (or $\Delta^{\ast}$-) ideal $P$ of the family  $(\eta_{i})_{i\in I}$ over $R$ is
a prime $\Delta$- (or $\Delta^{\ast}$- if we consider differences or inversive differences) ideal of the ring $R\{(y_{i})_{i\in
I}\}$ (respectively, of the ring of $\Delta^{\ast}$-polynomials $R\{(y_{i})_{i\in I}\}^{\ast}$). Therefore, $R\langle(\eta_{i})_{i\in I}\rangle$ can be treated as the quotient $\Delta$-field of the $\Delta$-ring $R\{(y_{i})_{i\in I}\}/P$. (In the case of inversive difference rings, the  $\Delta^{\ast}$-field $R\langle(\eta_{i})_{i\in I}\rangle^{\ast}$ can be considered as a quotient $\Delta^{\ast}$-field of the $\Delta^{\ast}$-ring
$R\{(y_{i})_{i\in I}\}^{\ast}/P$.)

\medskip

Let $K$ be a $\Delta$- (or $\Delta^{\ast}$-) field and $s$ a positive integer. By an {\em $s$-tuple over} $K$ we mean an $s$-dimensional vector $a = (a_{1},\dots, a_{s})$ whose coordinates belong to some $\Delta$- (respectively, $\Delta^{\ast}$-) overfield of $K$.

\begin{definition} Let $K$ be a $\Delta$- (or $\Delta^{\ast}$-) field and let $R$ be the algebra of $\Delta$- (respectively, $\Delta^{\ast}$-) polynomials in finitely many $\Delta$- (respectively, $\Delta^{\ast}$-) indeterminates $y_{1},\dots, y_{s}$ over $K$. Furthermore, let
$\Phi = \{f_{j} | j\in J\}$ be a set of $\Delta$- (respectively, $\Delta^{\ast}$-) polynomials in $R$. An $s$-tuple $\eta =
(\eta_{1},\dots, \eta_{s})$ over $K$ is said to be a {\bf solution} of the set $\Phi$ or a solution of the system of algebraic $\Delta$- (respectively, $\Delta^{\ast}$-) equations  $f_{j}(y_{1},\dots, y_{s}) = 0$ ($j\in J$) if $\Phi$ is contained in the kernel of the substitution of $(\eta_{1},\dots, \eta_{s})$ for $(y_{1},\dots, y_{s})$. In this case we also say that $\eta$ annuls $\Phi$.
\end{definition}
A system of algebraic difference equations $\Phi$ is called {\em prime} if the $\Delta$-ideal (or $\Delta^{\ast}$-ideal in the case of a system of difference or inversive difference equations) generated by $\Phi$ in the ring of $\Delta$ (or $\Delta^{\ast}$- if we deal with inversive difference equations) polynomials is prime.

As we have seen, if one fixes an $s$-tuple $\eta = (\eta_{1},\dots, \eta_{s})$ over a $\Delta$- (or $\Delta^{\ast}$-) field $K$, then all $\Delta$- (respectively, $\Delta^{\ast}$-) polynomials of the ring $K\{y_{1},\dots, y_{s}\}$ (respectively, $K\{y_{1},\dots, y_{s}\}^{\ast}$), for which $\eta$ is a solution, form a prime $\Delta$- (respectively, $\Delta^{\ast}$-) ideal, the defining $\Delta$- (respectively, $\Delta^{\ast}$-)  ideal of $\eta$. If $\Phi$ is a subset of $K\{y_{1},\dots, y_{s}\}$ (respectively, $K\{y_{1},\dots, y_{s}\}^{\ast}$), then an $s$-tuple
$\eta = (\eta_{1},\dots, \eta_{s})$ over $K$ is called a {\em generic zero} of $\Phi$ if for any  $\Delta$- (respectively, $\Delta^{\ast}$-) polynomial $f$, the inclusion $f\in \Phi$ holds if and only if $f(\eta_{1},\dots, \eta_{s})=0$.

Two $s$-tuples $\eta = (\eta_{1},\dots, \eta_{s})$ and $\zeta = (\zeta_{1},\dots, \zeta_{s})$ over a $\Delta$- (or $\Delta^{\ast}$-) field $K$ are called {\em equivalent} over $K$ if there is a $\Delta$-homomorphism $K\langle \eta_{1},\dots, \eta_{s}\rangle
\rightarrow K\langle \zeta_{1},\dots, \zeta_{s}\rangle$ (respectively, $K\langle \eta_{1},\dots, \eta_{s}\rangle^{\ast}
\rightarrow K\langle \zeta_{1},\dots, \zeta_{s}\rangle^{\ast}$) that maps each $\eta_{i}$ onto $\zeta_{i}$ and leaves the field $K$ fixed.

\begin{proposition}{\em (see {\em (\cite[Chapter 2, Theorem VII]{Cohn},
\cite[Propositions 3.2.6,  3.3.7]{KLMP})}\em}.  Let $R$ denote the algebra of $\Delta$- (or $\Delta^{\ast}$-) polynomials in $s$ $\Delta$- (respectively, $\Delta^{\ast}$-) indeterminates $y_{1},\dots, y_{s}$ over a $\Delta$- (respectively, $\Delta^{\ast}$-) field $K$.

{\em (i)}\, A set $\Phi \subsetneqq R$ has a generic zero if and only if $\Phi$ is a prime $\Delta$- (or $\Delta^{\ast}$-, if we consider differences or inversive differences) ideal of $R$. If $(\eta_{1},\dots, \eta_{s})$ is a generic zero of $\Phi$, then $K\langle \eta_{1},\dots, \eta_{s}\rangle$ is $\Delta$-isomorphic to the  $\Delta$- (respectively, $\Delta^{\ast}$-) quotient field of $R/\Phi$.

{\em (ii)}\, Any $s$-tuple over $K$ is a generic zero of some prime $\Delta$- (or $\Delta^{\ast}$-, if we deal with difference or inversive difference polynomials) ideal of $R$.  If two $s$-tuples over $K$ are generic zeros of the same prime $\Delta$- (or $\Delta^{\ast}$-) ideal of $R$, then these $s$-tuples are equivalent.
\end{proposition}

{\bf 2.3.\, Ring of differential, difference, and inversive difference operators. Differential, difference, and inversive difference modules.} \,

\medskip

Let $R$ be a differential or difference ring with a basic set $\Delta = \{\delta_{1},\dots, \delta_{m}\}$ and let $\Theta$ be the free commutative
semigroup generated by $\Delta$. If  $\theta = \delta_{1}^{k_{1}}\dots \delta_{m}^{k_{m}}\in \Theta$ ($k_{1},\dots, k_{m}\in {\mathbb N}$), then the number $ord\,\theta = \sum_{\nu = 1}^{m}k_{\nu}$ is called the {\em order} of $\theta$. Furthermore, for any $r\in {\mathbb N}$, the set $\{\theta\in\Theta\,|\,ord\,\theta \leq r\}$ is denoted by $\Theta(r)$.

\begin{definition} An expression of the form $\sum_{\theta\in \Theta}a_{\theta}\theta$, where $a_{\theta}\in R$ for any
$\theta\in \Theta$ and only finitely many elements $a_{\theta}$ are different from $0$, is called a $\Delta$-operator over the ring $R$. (If $\Delta$ is the set of mutually commuting derivations, then a $\Delta$-operator is also called a {\em differential operator}; if $\Delta$ consists of mutually commuting injective endomorphisms, a $\Delta$-operator is called a {\em difference operator}). Two $\Delta$-operators $\sum_{\theta\in
\Theta}a_{\theta}\theta$ and $\sum_{\theta\in\Theta}b_{\theta}\theta$ are considered to be equal if and only if $a_{\theta} = b_{\theta}$ for all $\theta\in \Theta$.
\end{definition}

The set of all $\Delta$-operators over a $\Delta$-ring $R$ can be equipped with a ring structure if we set \, $\sum_{\theta\in
\Theta}a_{\theta}\theta + \sum_{\theta\in \Theta}b_{\theta}\theta = \sum_{\theta\in \Theta}(a_{\theta} + b_{\theta})\theta$, $a\sum_{\theta\in \Theta}a_{\theta}\theta = \sum_{\theta\in \Theta}(aa_{\theta})\theta$, $(\sum_{\theta\in \Theta}a_{\theta}\theta)\theta_{1} = \sum_{\theta\in \Theta}a_{\theta}(\theta\theta_{1})$, $\delta a = a\delta + \delta(a)$ (respectively, $\delta a = \delta(a)\delta$ if $R$ is a difference ring and $\Delta$ is the basic set of endomorphisms of $R$) for any $\Delta$-operators $\sum_{\theta\in\Theta}a_{\theta}\theta,\,$ $\sum_{\theta\in \Theta}b_{\theta}\theta$ and for any $a\in R$, $\delta\in \Delta$, and extend the multiplication by distributivity. The ring obtained in this way is called {\em the ring of $\Delta$-operators}  over $R$; it will be denoted by ${\cal{D}}$. (If $\Delta$ is a set of derivations, ${\cal{D}}$ is also  said to be the ring of differential operators over the differential ring $R$; if $\Delta$ is a set of endomorphisms, ${\cal{D}}$ is called the ring of difference operators over $R$.)

\medskip

The order of a nonzero $\Delta$-operator  $A = \sum_{\theta\in \Theta}a_{\theta}\theta \in {\cal{D}}$ is defined as the number $ord\,A = \max\{ord\,\theta\,|\,a_{\theta}\neq 0\}$. We also set $ord\,0 = -\infty$.

Let ${\cal{D}}_{r} = \{A\in {\cal{D}} | ord\,A \leq r\}$ for any $r\in{\mathbb N}$ and let ${\cal{D}}_{r} = 0$ for any $r\in {\mathbb Z}, r < 0$. Then the ring $\cal{D}$ can be treated as a filtered ring with the ascending filtration $({\cal{D}}_{r})_{r\in{\mathbb Z}}$. Below, while considering $\cal{D}$ as a filtered ring, we always mean this filtration.

\begin{definition}
Let $R$ be a $\Delta$-ring and $\cal{D}$ the ring of $\Delta$-operators over $R$. Then a left $\cal{D}$-module is called a $\Delta$-$R$-module.  {\em (If   $\Delta$ is a set of derivations, we also use the term} differential $R$-module; {\em if  $\Delta$ is a set of endomorphisms, we use the term} difference $R$-module{\em )}. In other words, an $R$-module $M$ is a $\Delta$-$R$-module if the elements of $\Delta$ act on $M$ in such a way that $\delta(x+y) = \delta(x)+\delta(y)$, $\delta(\delta' x) = \delta'(\delta x)$, and $\delta(ax)=\delta(a)x + a\delta(x)$ (if $\Delta$ consists of derivations, so $\Delta$- means ''differential'') or $\delta(ax)=\delta(a)\delta(x)$ (if $\Delta$ consists of endomorphisms, so $\Delta$- means ''difference'') for any $x, y\in M; \delta,\, \delta'\in\Delta$; $a\in R$. If $R$ is a $\Delta$-field, then a $\Delta$-$R$-module $M$ is also called a  vector $\Delta$-$R$-space.
\end{definition}

If $R$ is an inversive difference ring with basic set $\Delta = \{\delta_{1},\dots, \delta_{m}\}$ and $\Gamma$ is the free commutative
group generated by $\Delta$, then the order of an element $\gamma = \delta_{1}^{k_{1}}\dots \delta_{m}^{k_{m}}\in \Gamma$ ($k_{1},\dots, k_{m}\in {\mathbb Z}$) is defined as $ord\,\gamma = \sum_{\nu = 1}^{m}|k_{\nu}|$. Also, for any $r\in {\mathbb N}$, we set $\Gamma(r) = \{\gamma\in\Gamma\,|\,ord\,\gamma \leq r\}$. In this case, by a $\Delta^{\ast}$-operator we mean an expression of the form $\sum_{\gamma\in\Gamma}a_{\gamma}\gamma$, where $a_{\gamma}\in R$ for any $\gamma\in\Gamma$ and only finitely many elements $a_{\gamma}$ are different from $0$. As in the case of $\Delta$-operators, two $\Delta^{\ast}$-operators are considered to be equal if and only if for any $\gamma\in \Gamma$, the coefficients of $\gamma$ in these operators are the same.

Clearly, a $\Delta^{\ast}$-ring $R$ can be also treated as a $\Delta$-ring, so every $\Delta$-operator over $R$ can be also considered as a $\Delta^{\ast}$-operator. The set of all $\Delta^{\ast}$-operators over the ring $R$ can be naturally considered as a ring extension of the ring ${\cal{D}}$ of $\Delta$-operators over $R$ where the operation are defined in the same way as they are defined in ${\cal{D}}$ with additional rules $\delta^{-1}a = \delta^{-1}(a)\delta^{-1}$ and $\delta\delta^{-1} = \delta^{-1}\delta = 1$ ($a\in R$, $\delta\in\Delta$) extended by distributivity. The resulting ring will be denoted by ${\cal{D}}^{\ast}$; it is called the ring of $\Delta^{\ast}$-operators (or {\em inversive difference operators}) over $R$. A left ${\cal{D}}^{\ast}$-module is called a $\Delta^{\ast}$-$R$-module (or an inversive difference $R$-module). Such a module is actually a $\Delta$-$R$-module $M$ with an additional action of the elements of the form $\delta^{-1}$ ($\delta\in\Delta$) such that $\delta(\delta^{-1}(x)) = \delta^{-1}(\delta(x))$ for any $x\in M$ (the other rules are the same as in Definition II.5 except for that the elements $\delta$ and $\delta'$ are taken from the set $\Delta^{\ast}$ rather than from $\Delta$).

The order of a nonzero $\Delta^{\ast}$-operator  $A = \sum_{\gamma\in \Gamma}a_{\gamma}\gamma \in {\cal{D}}^{\ast}$ is defined as the number $ord\,A = \max\{ord\,\gamma\,|\,a_{\gamma}\neq 0\}$, and we also set $ord\,0 = -\infty$. The ring ${\cal{D}}^{\ast}$ will be treated as a  filtered ring with the ascending filtration $({\cal{D}}^{\ast}_{r})_{r\in{\mathbb Z}}$ where ${\cal{D}}^{\ast}_{r} = \{A\in {\cal{D}}^{\ast} | ord\,A \leq r\}$ for any $r\in{\mathbb N}$ and ${\cal{D}}^{\ast}_{r} = 0$ for any $r\in {\mathbb Z}, r < 0$.

\section{Differential and difference dimension polynomials}

In this section we present main theorems on dimension polynomials of differential and difference modules and field extensions. Then we show how one can determine the strength of a system of partial differential or difference equations by computing the corresponding dimension polynomial.

\medskip

With the above notation, let $R$ be a $\Delta$- (respectively, $\Delta^{\ast}$-) ring. We say that a $\Delta$-$R$-module (respectively, $\Delta^{\ast}$-$R$-module) $M$ is finitely generated, if it is finitely generated as a left $\cal{D}$- (respectively, ${\cal{D}}^{\ast}$-) module. By a filtered $\Delta$- (respectively, $\Delta^{\ast}$-) module we always mean a left $\cal{D}$- (respectively, ${\cal{D}}^{\ast}$-) module $M$ equipped with an exhaustive and separated filtration. Thus, a filtration of $M$ is an ascending chain $(M_{r})_{r\in{\mathbb Z}}$ of $R$-submodules of $M$ such that ${\cal{D}}_{r}M_{s}\subseteq M_{r+s}$ (respectively, ${\cal{D}}^{\ast}_{r}M_{s}\subseteq M_{r+s}$) for all $r, s\in {\mathbb Z}$, $M_{r} = 0$ for all sufficiently small $r\in {\mathbb Z}$, and $\bigcup_{r\in {\mathbb Z}}M_{r} = M$. A filtration $(M_{r})_{r\in{\mathbb Z}}$ of $M$ is called {\em excellent} if all $R$-modules $M_{r}$ (${r\in{\mathbb Z}}$) are finitely generated and there exists $r_{0}\in{\mathbb Z}$ such that $M_{r} = {\cal{D}}_{r-r_{0}}M_{r_{0}}$ (respectively, $M_{r} = {\cal{D}}^{\ast}_{r-r_{0}}M_{r_{0}}$) for any ${r\in{\mathbb Z}}, r \geq r_{0}$. Note that if $R$ is a $\Delta$-field and $M$ is a finitely generated $\Delta$-$R$-module, $$M = \D\sum_{i=1}^{s}{\cal{D}}f_{i}$$ for some $f_{1},\dots, f_{s}\in M$, then $$\left(M_{r} = \D\sum_{i=1}^{s}{\cal{D}}_{r}f_{i}\right)_{r\in{\mathbb Z}}$$ is an excellent filtration of $M$, and a similar remark can be made about a finitely generated $\Delta^{\ast}$-$R$-module.

The following result combines theorems on dimension polynomials of differential and difference modules obtained in \cite{Johnson1} and \cite{Levin1} (see also \cite[Theorems 5.1.11, 6.2.5 and Propositions 5.2.12, 6.2.17]{KLMP}).

\begin{theorem}
Let $R$ be a $\Delta$-field whose basic set consists of $m$ operators (deri\-vations or injective endomorphisms). Let $\cal{D}$ be the ring of $\Delta$-operators over $R$, and let $(M_{r})_{r\in{\mathbb Z}}$ be an excellent filtration of a $\Delta$-$R$-module $M$. Then there exists a polynomial $\psi(t)\in {\mathbb Q}[t]$ with the following properties.

{\em (i)}\, $\psi(r) = \dim_{R}(M_{r})$ for all sufficiently large $r\in {\mathbb Z}$, that is, there exists $r_{0}\in {\mathbb Z}$ such that the last equality holds for all integers $r\geq r_{0}$.  {\em (as usual, $\dim_{R}(M_{r})$ denotes the dimension of the vector $R$-space $M_{r}$)}.

{\em (ii)}\, $deg\,\psi(t)\leq m$ and the polynomial $\psi(t)$ can be written as
 $$\psi(t) = \sum_{i=0}^{m}c_{i}\binom{t+i}{i}$$ where
$c_{0}, c_{1},\dots, c_{m}\in {\mathbb Z}$ and $$\D\binom{t+i}{i} = \frac{(t+i)(t+i-1)\ldots (t+1)}{i!}$$ {\em (this polynomial takes integer values for all sufficiently large integer values of $t$).}

{\em (iii)}\, The integers $d=deg\,\psi(t),\, c_{m}$\, and $c_{d}$ (if $d < m$) do not depend on the choice of the excellent filtration of $M$. Furthermore, $c_{m}$ is equal to the maximal number of elements of $M$ linearly independent over the ring $\cal{D}$.
\end{theorem}

The polynomial $\psi(t)$ whose existence is established by Theorem III.1 is called the {\em $\Delta$-dimension polynomial} (differential or difference dimension polynomial depending on the nature of the set $\Delta$) of the $\Delta$-$R$-module $M$ associated with the excellent filtration $(M_{r})_{r\in{\mathbb Z}}$. The integers $d$, $c_{m}$, and $c_{d}$ are called the {\em $\Delta$-type,  $\Delta$-dimension}, and {\em typical $\Delta$-dimension} of $M$, respectively. A number of results on differential and difference dimension polynomials, as well as some methods of their computation, can be found in \cite[Chapters 5 - 9]{KLMP}.

The following is an analog of Theorem III.1 for inversive difference modules (see \cite[Theorem 6.3.3 and Proposition 6.3.15]{KLMP} or \cite[Theorems 3.5.2, 3.5.8]{Levin7}).

\begin{theorem}
Let $R$ be a $\Delta^{\ast}$-field  whose basic set $\Delta$ consists of $m$ automorphisms of $R$, and let $(M_{r})_{r\in{\mathbb Z}}$ be an excellent filtration of a $\Delta^{\ast}$-$R$-module $M$. Then there exists a polynomial $\chi(t)$ in one variable $t$ with rational coefficients such that

\smallskip

{\em (i)}\, $\chi(r) = \dim_{R}(M_{r})$ for all sufficiently large $r\in{\mathbb Z}$\,;

\smallskip

{\em (ii)}\, $deg\,\chi(t)\leq m$ and the polynomial $\chi(t)$ can be represented in the form
$$\chi(t) = {\frac{2^{m}a}{m!}}t^{m} + o(t^{m})$$ where $a\in {\mathbb Z}$ and $o(t^{n})$ is a polynomial in ${\mathbb Q}[t]$ of degree less than $m$.

\smallskip

{\em (iii)}\, The integers $a$, $d = deg\,\chi(t)$ and the coefficient of $t^{d}$ in the polynomial $\chi(t)$ do not depend on the choice of the excellent filtration of $M$. Furthermore, $a$ is equal to the maximal number of elements of $M$ linearly independent over the ring ${\cal{D}}^{\ast}$.

\end{theorem}

The polynomial $\chi(t)$ is called the {\em $\Delta^{\ast}$-dimension polynomial} of the $\Delta^{\ast}$-$R$-module $M$ associated with the excellent filtration $(M_{r})_{r\in{\mathbb Z}}$.

\smallskip

The next result combines Kolchin's theorem on differential dimension polynomial
\cite[Chapter II, Theorem 6]{Kolchin} and the corresponding result for difference field extensions proved in \cite{Levin1}.

\begin{theorem}
Let $K$ be a $\Delta$-field whose basic set consists of $m$ operators (derivations or endomorphisms).  Let $L = K\langle \eta_{1},\dots, \eta_{s}\rangle$ be a $\Delta$-field extension of $K$ generated by a finite family $\eta = \{\eta_{1},\dots, \eta_{s}\}$. Then there exists a polynomial $\phi_{\eta | K}(t)\in {\mathbb Q}[t]$ with the following properties.

\smallskip

{\em (i)\,} $\phi_{\eta | K}(r) = trdeg_{K}K(\{\theta\eta_{j}\,|\, \theta \in \Theta(r), 1\leq j\leq s\})$ {\em for all sufficiently large}
$r\in {\mathbb N}$.

{\em (ii)\,}  $deg\,\phi_{\eta | K}(t)\leq n$ and the polynomial $\phi_{\eta | K}(t)$ can be written as $$\phi_{\eta | K}(t) =
\sum_{i=0}^{m}a_{i}\binom{t+i}{i}$$ where $a_{0},\dots, a_{m}\in {\mathbb Z}$.

\smallskip

 {\em (iii)\,} The integers $a_{m}, d = deg\,\phi_{\eta | K}(t)$ and $a_{d}$ are invariants of the polynomial $\phi_{\eta | K}(t)$, that is, they do  not depend on the choice of a system of $\sigma$-generators $\eta$. Furthermore, $a_{m} = \Delta$-$trdeg_{K}L$ where $\Delta$-$trdeg_{K}L$ denotes the  $\Delta$-transcendence degree of $L$ over $K$, that is, the maximal number of elements $\xi_{1},\dots,\xi_{k}\in L$ such that the family
$\{\theta\xi_{i} | \theta\in\Theta, 1\leq i\leq k\}$ is algebraically independent over $K$.
\end{theorem}

The polynomial $\phi_{\eta | K}(t)$ whose existence is established by Theorem III.3 is called the
 $\Delta$- ({\em differential} or {\em difference} depending on the nature of $\Delta$) {\em dimension polynomial of the $\Delta$-field
extension $L$ of $K$ associated with the system of $\Delta$-generators $\eta$}. The integers $d = deg\,\phi_{\eta\,|\,K}(t)$ and $a_{d}$ are called, respectively, the $\Delta$-{\em type} and {\em typical $\Delta$-transcendence degree} of $L$ over $K$. These invariants of $\phi_{\eta | K}(t)$ are denoted by $\Delta$-$type_{K}L$ and $\Delta$-$t.trdeg_{K}L$, respectively.

\medskip

Notice that if the elements $\eta_{1},\dots, \eta_{s}$ are $\Delta$-algebraically independent over $K$ (that is, the set $\{\theta\eta_{i}\,|\,\theta\in \Theta,\, 1\leq i\leq s\}$ is algebraically independent over $K$) and $\phi_{\eta | K}(t)$ is the corresponding
difference dimension polynomial of $L/K$ (we use the notation of the last theorem), then $$\phi_{\eta | K}(r) = trdeg_{K}K(\{\tau\eta_{j}\,|\,\tau\in T,\, 1\leq j\leq s\}) = s\cdot\,Card\,\Theta(r) = s\D\binom{r+m}{m}$$ for all sufficiently large $r\in {\mathbb N}$ ($Card\,\Theta(r)$ is the number of solutions $(k_{1},\dots, k_{m})\in {\mathbb N}^{m}$ of the inequality $k_{1}+\dots + k_{m}\leq r$; it is well known (see, for example, \cite[Proposition 2.1.9]{KLMP}) that this number is ${\binom{r+m}{m}}$). Therefore, in this case $$\phi_{\eta | K}(t) = s\D{\binom{t+m}{m}}.$$

\medskip

The following theorem shows the existence of a dimension polynomial of a finitely generated inversive difference field extension.

\begin{theorem} Let $K$ be a $\Delta^{\ast}$-field  whose basic set $\Delta$ consists of $m$ automorphisms of $K$. As before, let $\Gamma$ be the free commutative group generated by $\Delta$, and for any $r\in {\mathbb N}$, let $\Gamma(r) = \{\gamma\in \Gamma\,|\,ord\,\gamma \leq r\}$. Furthermore, let $L
= K\langle \eta_{1},\dots, \eta_{s}\rangle^{\ast}$ be a $\Delta^{\ast}$-field extension of $K$ generated by a finite family $\eta = \{\eta_{1},\dots, \eta_{s}\}$. Then there exists a polynomial $\psi_{\eta | K}(t)\in {\mathbb Q}[t]$ with the following properties.

{\em (i)}\,  $\psi_{\eta | K}(r) = trdeg_{K}K(\{\gamma \eta_{j}\,|\,\gamma \in \Gamma(r), 1\leq j\leq s\})$  for all sufficiently
large $r\in {\mathbb N}$.

{\em (ii)}\, $deg\,\psi_{\eta | K}(t)\leq m$ and the polynomial $\psi_{\eta | K}(t)$ can be written as $$\psi_{\eta | K}(t) = {\frac{2^{m}a}{m!}}t^{m} + o(t^{m})$$ where $a\in {\mathbb Z}$ and $o(t^{m})$ is a polynomial in ${\mathbb Q}[t]$ of degree less than $m$.

\smallskip

{\em (iii)}\, The integers $a$, $d = deg\,\psi_{\eta | K}(t)$ and the coefficient of $t^{d}$ in the polynomial $\psi_{\eta | K}(t)$ do not depend on the choice of a system of generators $\eta$. Furthermore, $a = \Delta$-$trdeg_{K}L$.

\smallskip

{\em (iv)}\, If $\eta_{1},\dots, \eta_{s}$ are $\Delta$-algebraically independent over $K$, then $$\psi_{\eta | K}(t) = s\D\sum_{k=0}^{m}(-1)^{m-k}2^{k}\binom{m}{k}{\binom{t+k}{k}}\,.$$

\end{theorem}

\medskip

Let us consider a prime system of algebraic $\Delta$- (differential or  difference) or $\Delta^{\ast}$-(inversive difference)
equations
\begin{equation}
A_{i}(y_{1},\dots, y_{s}) = 0\hspace{0.5in}(i=1,\dots, p)
\end{equation}
where $A_{i}(y_{1},\dots, y_{s})$ are  $\Delta$- (or $\Delta^{\ast}$- ) polynomials in the ring $R = K\{y_{1},\dots, y_{s}\}$ (respectively, in $R= K\{y_{1},\dots, y_{s}\}^{\ast}$) and let $P$ be a prime $\Delta$-ideal (respectively a prime $\Delta^{\ast}$-ideal if we consider the difference or inversive difference case) of $R$ generated by the right-hand sides of system (1). Furthermore, let $\eta_{i}$ be the canonical image of $y_{i}$ in the factor ring $R/P$ ($1\leq i\leq s$). It is easy to see that for every $r\in {\mathbb N}$ the intersection $P\cap R_{r}$ is a prime ideal of the ring $R_{r}$ and the quotient fields of the rings $R_{r}/P\cap R_{r}$ and $K[\{\theta(\eta_{j})\,|\,\theta\in\Theta(r), \,1\leq j\leq s\}]$ (respectively, $K[\{\gamma(\eta_{j})\,|\,\gamma\in \Gamma(r), \,1\leq j\leq s\}]$) are isomorphic.  Considering the case of algebraic differential or difference equations we can apply Theorem III.3 and obtain that there exists a polynomial $\phi_{P}(t)$ in one variable $t$ with rational coefficients such that $$\phi_{P}(t) = trdeg_{K}K(\{\theta(\eta_{j})\,|\,\theta\in\Theta(r), \,1\leq j\leq s\}) = trdeg_{K}(R_{r}/P\cap R_{r})$$ for all sufficiently large
$r\in {\mathbb Z}$,\, $deg\,\phi_{P}(t)\leq m$ and the polynomial $\phi_{P}(t)$ can be written as $$\phi_{P}(t) =
\sum_{i=0}^{m}a_{i}\binom{t+i}{i}$$ where $a_{0},\dots, a_{m}\in {\mathbb Z}$ and $a_{m} = \Delta$-$trdeg_{K}K(\{\theta(\eta_{j})\,|\,\theta\in\Theta(r), \,1\leq j\leq s\})$.

In the case of a system of difference equations (including the case when such a system involves negative degrees of basic translations, which act as automorphisms), one can apply Theorem III.4 that shows the existence of a polynomial $\psi_{P}(t)\in{\mathbb Q}[t]$ such that $$\psi_{P}(r) = trdeg_{K}K(\{\gamma(\eta_{j})\,|\,\gamma\in \Gamma(r), \,1\leq j\leq s\}) = trdeg_{K}(R_{r}/P\cap R_{r})$$ for all sufficiently large $r\in {\mathbb Z}$,\, $deg\,\psi(t)\leq m$ and the polynomial $\psi_{P}(t)$ can be written as $$\psi_{P}(t) = {\frac{2^{n}a_{P}}{n!}}t^{n} + o(t^{n})$$ where $a_{P} = \sigma$-$trdeg_{K}(R/P)$.

\bigskip

With the above notation, the numerical polynomial $\phi_{P}(t)$ (respectively, $\psi_{P}(t)$) is called a {\em differential} (respectively, {\em difference}) {\em dimension polynomial} of system (1). It is also said to be a {\em $\Delta$- (respectively, $\Delta^\ast$-) dimension polynomial} of the system.

\medskip

Taking into account Einstein's approach described in the Introduction, one can say that for all sufficiently large $r$, the value $\phi_{P}(r)$ of the differential dimension polynomial $\phi_{P}(t)$ of a system of algebraic differential equations is the number of Taylor coefficients of order $\leq r$ of an analytic solution that can be chosen arbitrarily. (These Taylor coefficients are the values of derivatives of order $\leq r$  of the solution computed at the point in whose neighborhood we consider its expansion. The dependence of coefficients is understood as their algebraic dependence over the field of coefficients of the system.) Thus, $\phi_{P}(t)$ can be viewed as a measure of strength of system (1), so
the problem of computation of differential dimension polynomials is important not only for the study of differential algebraic structures, but also for the study of equations of mathematical physics.

Considering a system of equations in finite differences over a field of functions in several real variables, one can use Einstein's approach to define the concept of {\em strength} of such a system as follows (cf. Einstein's description of the strength of a system of PDEs presented in the Introduction).  Let

\begin{equation}
A_{i}(f_{1},\dots, f_{s}) = 0\hspace{0.3in}(i=1,\dots, p)
\end{equation}
be a system of equations in finite differences with respect to $s$ unknown grid functions $f_{1},\dots, f_{s}$ in $n$ real variables $x_{1},\dots, x_{n}$ with coefficients in some functional field $K$. We also assume that the difference grid, whose nodes form the domain of considered functions, has equal cells of dimension $h_{1}\times\dots\times h_{n}$ \,($h_{1},\dots, h_{n}\in {\mathbb R}$) and fills the whole space ${\mathbb R}^{n}$. As an example, one can consider a field $K$ consisting of a zero function and fractions of the form $u/v$ where $u$ and $v$ are grid functions defined almost everywhere and vanishing at at most finitely many nodes. (As usual, we say that a grid function is defined almost everywhere if there are at most
finitely many nodes where it is not defined.)

Let us fix some node $\mathcal{P}$ and  say that {\em a node $\mathcal{Q}$ has order $i$} (with respect to $\mathcal{P}$) if the shortest path from $\mathcal{P}$ to $\mathcal{Q}$ along the edges of the grid consists of $i$ steps (by a step we mean a path from a node of the grid to a neighboring node along the edge between these two nodes). Say, the orders of the nodes in the two-dimensional case are as follows (a number near a node shows the order of this node).
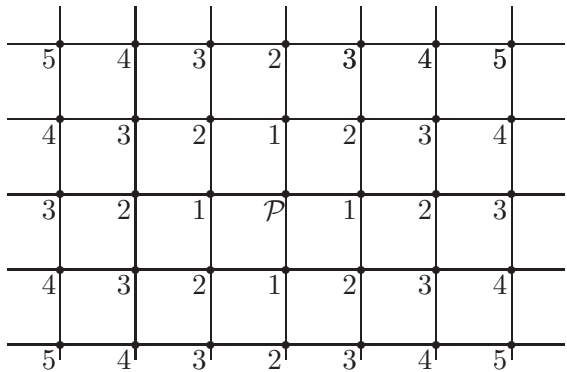
\begin{figure}[h]
\setlength{\unitlength}{1cm}
\begin{picture}(9,4.9)
\put(1.6,0.3){\line(1,0){0.7}} \put (2.3,0.30){\circle*{0.1}}
\put(2.3,0.3){\line(1,0){1}} \put (3.3,0.30){\circle*{0.1}} \put
(4.3,0.30){\circle*{0.1}} \put (6.3,0.30){\circle*{0.1}} \put
(5.3,0.30){\circle*{0.1}} \put(3.3,0.3){\line(1,0){1}}
\put(4.3,0.3){\line(1,0){1}} \put(5.3,0.3){\line(1,0){1}}
\put(6.3,0.3){\line(1,0){1}} \put (7.3,0.30){\circle*{0.1}}
\put(7.3,0.3){\line(1,0){1}} \put (8.3,0.30){\circle*{0.1}}
\put(8.3,0.3){\line(1,0){0.7}} \put(2.3,0.3){\line(0,-1){0.2}}
\put(3.3,0.3){\line(0,-1){0.2}} \put(4.3,0.3){\line(0,-1){0.2}}
\put(5.3,0.3){\line(0,-1){0.2}} \put(6.3,0.3){\line(0,-1){0.2}}
\put(7.3,0.3){\line(0,-1){0.2}} \put(8.3,0.3){\line(0,-1){0.2}}
\put(2.3,0.3){\line(0,1){1}} \put(3.3,0.3){\line(0,1){1}}
\put(4.3,0.3){\line(0,1){1}} \put(5.3,0.3){\line(0,1){1}}
\put(6.3,0.3){\line(0,1){1}} \put(7.3,0.3){\line(0,1){1}}
\put(8.3,0.3){\line(0,1){1}} \put(7.3,1.30){\circle*{0.1}}
\put(2.3,1.30){\circle*{0.1}} \put(3.3,1.30){\circle*{0.1}}
\put(4.3,1.30){\circle*{0.1}} \put(5.3,1.30){\circle*{0.1}}
\put(6.3,1.30){\circle*{0.1}} \put(8.3,1.30){\circle*{0.1}}
\put(1.6,1.3){\line(1,0){0.7}} \put(2.3,1.3){\line(1,0){1}}
\put(3.3,1.3){\line(1,0){1}} \put(4.3,1.3){\line(1,0){1}}
\put(5.3,1.3){\line(1,0){1}} \put(6.3,1.3){\line(1,0){1}}
\put(7.3,1.3){\line(1,0){1}} \put(8.3,1.3){\line(1,0){0.7}}

\put(2.3,2.3){\circle*{0.1}} \put(3.3,2.3){\circle*{0.1}}
\put(4.3,2.3){\circle*{0.1}} \put(5.3,2.3){\circle*{0.1}}
\put(6.3,2.3){\circle*{0.1}} \put(7.3,2.3){\circle*{0.1}}
\put(8.3,2.3){\circle*{0.1}}

\put(2.3,3.3){\circle*{0.1}} \put(3.3,3.3){\circle*{0.1}}
\put(4.3,3.3){\circle*{0.1}} \put(5.3,3.3){\circle*{0.1}}
\put(6.3,3.3){\circle*{0.1}} \put(7.3,3.3){\circle*{0.1}}
\put(8.3,3.3){\circle*{0.1}}

\put(2.3,4.3){\circle*{0.1}} \put(3.3,4.3){\circle*{0.1}}
\put(4.3,4.3){\circle*{0.1}} \put(5.3,4.3){\circle*{0.1}}
\put(6.3,4.3){\circle*{0.1}} \put(7.3,4.3){\circle*{0.1}}
\put(8.3,4.3){\circle*{0.1}}

\put(1.6,2.3){\line(1,0){0.7}} \put(2.3,2.3){\line(1,0){1}}
\put(3.3,2.3){\line(1,0){1}} \put(4.3,2.3){\line(1,0){1}}
\put(5.3,2.3){\line(1,0){1}} \put(6.3,2.3){\line(1,0){1}}
\put(7.3,2.3){\line(1,0){1}} \put(8.3,2.3){\line(1,0){0.7}}

\put(1.6,3.3){\line(1,0){0.7}} \put(2.3,3.3){\line(1,0){1}}
\put(3.3,3.3){\line(1,0){1}} \put(4.3,3.3){\line(1,0){1}}
\put(5.3,3.3){\line(1,0){1}} \put(6.3,3.3){\line(1,0){1}}
\put(7.3,3.3){\line(1,0){1}} \put(8.3,3.3){\line(1,0){0.7}}

\put(1.6,4.3){\line(1,0){0.7}} \put(2.3,4.3){\line(1,0){1}}
\put(3.3,4.3){\line(1,0){1}} \put(4.3,4.3){\line(1,0){1}}
\put(5.3,4.3){\line(1,0){1}} \put(6.3,4.3){\line(1,0){1}}
\put(7.3,4.3){\line(1,0){1}} \put(8.3,4.3){\line(1,0){0.7}}

\put(2.3,1.3){\line(0,1){1}} \put(3.3,1.3){\line(0,1){1}}
\put(4.3,1.3){\line(0,1){1}} \put(5.3,1.3){\line(0,1){1}}
\put(6.3,1.3){\line(0,1){1}} \put(7.3,1.3){\line(0,1){1}}
\put(8.3,1.3){\line(0,1){1}}

\put(2.3,2.3){\line(0,1){1}} \put(3.3,2.3){\line(0,1){1}}
\put(4.3,2.3){\line(0,1){1}} \put(5.3,2.3){\line(0,1){1}}
\put(6.3,2.3){\line(0,1){1}} \put(7.3,2.3){\line(0,1){1}}
\put(8.3,2.3){\line(0,1){1}}

\put(2.3,3.3){\line(0,1){1}} \put(3.3,3.3){\line(0,1){1}}
\put(4.3,3.3){\line(0,1){1}} \put(5.3,3.3){\line(0,1){1}}
\put(6.3,3.3){\line(0,1){1}} \put(7.3,3.3){\line(0,1){1}}
\put(8.3,3.3){\line(0,1){1}}

\put(2.3,4.3){\line(0,1){0.5}} \put(3.3,4.3){\line(0,1){0.5}}
\put(4.3,4.3){\line(0,1){0.5}} \put(5.3,4.3){\line(0,1){0.5}}
\put(6.3,4.3){\line(0,1){0.5}} \put(7.3,4.3){\line(0,1){0.5}}
\put(8.3,4.3){\line(0,1){0.5}}

\put(5.15,2.1){\makebox(0,0){$\mathcal{P}$}}
\put(4.15,2.1){\makebox(0,0){$1$}}
\put(6.15,2.1){\makebox(0,0){$1$}}
\put(3.15,2.1){\makebox(0,0){$2$}}
\put(7.15,2.1){\makebox(0,0){$2$}}
\put(2.15,2.1){\makebox(0,0){$3$}}
\put(8.15,2.1){\makebox(0,0){$3$}}

\put(5.15,1.1){\makebox(0,0){$1$}}
\put(5.15,3.1){\makebox(0,0){$1$}}
\put(4.15,1.1){\makebox(0,0){$2$}}
\put(6.15,1.1){\makebox(0,0){$2$}}
\put(3.15,1.1){\makebox(0,0){$3$}}
\put(7.15,1.1){\makebox(0,0){$3$}}
\put(2.15,1.1){\makebox(0,0){$4$}}
\put(8.15,1.1){\makebox(0,0){$4$}}

\put(5.15,4.1){\makebox(0,0){$2$}}
\put(5.15,0.1){\makebox(0,0){$2$}}
\put(4.15,0.1){\makebox(0,0){$3$}}
\put(6.15,4.1){\makebox(0,0){$3$}}
\put(3.15,0.1){\makebox(0,0){$4$}}
\put(7.15,4.1){\makebox(0,0){$4$}}
\put(2.15,0.1){\makebox(0,0){$5$}}
\put(8.15,4.1){\makebox(0,0){$5$}}

\put(6.15,0.1){\makebox(0,0){$3$}}
\put(7.15,0.1){\makebox(0,0){$4$}}
\put(8.15,0.1){\makebox(0,0){$5$}}

\put(6.15,4.1){\makebox(0,0){$3$}}
\put(7.15,4.1){\makebox(0,0){$4$}}
\put(8.15,4.1){\makebox(0,0){$5$}}

\put(6.15,3.1){\makebox(0,0){$2$}}
\put(7.15,3.1){\makebox(0,0){$3$}}
\put(8.15,3.1){\makebox(0,0){$4$}}

\put(4.15,4.1){\makebox(0,0){$3$}}
\put(3.15,4.1){\makebox(0,0){$4$}}
\put(2.15,4.1){\makebox(0,0){$5$}}

\put(4.15,3.1){\makebox(0,0){$2$}}
\put(3.15,3.1){\makebox(0,0){$3$}}
\put(2.15,3.1){\makebox(0,0){$4$}}
\end{picture}
\caption{2-dimensional grid}
\end{figure}
Let us consider the values of the unknown grid functions $f_{1},\dots, f_{s}$ at the nodes whose order does not exceed $r$\, ($r\in {\mathbb N}$). If $f_{1},\dots, f_{s}$ should not satisfy any system of equations (or any other condition), their values at nodes of any order can be chosen arbitrarily. Because of the system in finite differences (and equations obtained from the equations of the system by transformations of the form $f_{j}(x_{1},\dots, x_{s})\mapsto f_{j}(x_{1}+k_{1}h_{1},\dots, x_{s}+k_{n}h_{n})$ with $k_{1},\dots, k_{n}\in {\mathbb Z}$, $1\leq j\leq s$), the number of independent values of the functions $f_{1},\dots, f_{s}$ at the nodes of order $\leq r$ decreases. This number, which is a function of $r$, is considered as a ''measure of strength'' of the system in finite differences (in the sense of Einstein). We denote it by $S_{r}$.\\
\indent With the above conventions, suppose that the transformations $\alpha_{j}$ of the field of coefficients $K$ defined by $$\alpha_{j}f(x_{1},\dots, x_{n}) = f(x_{1},\dots,x_{j-1}, x_{j}+h_{j},\dots, x_{n})$$ ($1\leq j\leq n$) are automorphisms of this field. Then $K$
can be considered as an inversive difference field with the basic set $\sigma = \{\alpha_{1},\dots, \alpha_{n}\}$.  The replacement of the unknown functions $f_{i}$ by difference indeterminates $y_{i}$ ($i=1,\dots, s$) leads to a system of algebraic difference equations of the form (1). If this system is prime (e.g., we deal with a system of linear difference equations), then its difference dimension polynomial $\psi(t)$ expresses the strength $S_{r}$. Thus, this polynomial can be naturally viewed as the measure of Einstein's strength of a given system of equations in finite differences. In what follows, the $\Delta^{\ast}$-dimension polynomial $\psi(t)$ will be called the {\em difference dimension polynomial of the system}.

\bigskip

Methods of computation of $\Delta$- and $\Delta^{\ast}$- dimension polynomials of a system of algebraic partial differential or difference equations developed so far are based either on building of a characteristic set of the considered above associated $\Delta$- (or $\Delta^{\ast}$- ) ideal $P$ in  $K\{y_{1},\dots, y_{s}\}$  or on constructing a free resolution of the module of K\"ahler differentials associated with the extension $K\langle\eta_{1},\dots, \eta_{s}\rangle$ (or $K\langle\eta_{1},\dots, \eta_{s}\rangle^{\ast}$). The corresponding computations can be found, for example, in \cite[Chapter 9]{KLMP} and \cite[Chapter 7]{Levin7}. The main drawback of the mentioned approaches is the lack of efficient algorithms for constructing characteristic sets and serious restrictions on the systems to which one can apply the method of free resolutions. In the last case, a system of difference equations with inversive difference operators is supposed to be linear and symmetric, that is, whenever an equation involves a $\Delta^{\ast}$-operator  $\omega = a_{1}\delta_{1}^{k_{11}}\dots \delta_{m}^{k_{1m}} + \dots + a_{r}\delta_{1}^{k_{r1}}\dots \delta_{m}^{k_{rm}}$ ($a_{i}\in K$), which contains a term $a\delta_{1}^{l_{1}}\dots \delta_{m}^{l_{m}}$ $(a\in K,\,a\neq 0)$, then it also contains all terms of the form
$b\delta_{1}^{\pm l_{1}}\dots \delta_{m}^{\pm l_{m}}$ with nonzero coefficients $b\in K$ and all $2^{m}$ distinct combinations of signs before $l_{1},\dots, l_{m}$. In what follows we explain a method of computation of dimension polynomials (and therefore, the strength of a system of
algebraic partial differential or difference equations), which does not have these restrictions.

Implementations for computing Gr\"obner bases in modules of differential and difference operators are available, e.g., in the Mgfun package \cite{Mgfun} for Maple or in the Plural extension of Singular \cite{Singular}.

\section{Computation of the strength of a system of difference equations via Gr\"obner and generalized Gr\"obner basis techniques. Examples}

Let $K$ be a difference or inversive difference field of characteristic $0$ with basic set $\Delta=\{\delta_1,\ldots,\delta_m\}$. As we have seen, the ring of $\Delta$-operators over $K$ carries many properties of a polynomial ring in $m$ variables over $K$. In order to underline the relationship between the Gr\"obner basis method for $\Delta$-$K$-modules considered below and the classical Gr\"obner basis technique for polynomial ideals we will denote the ring of $\Delta$-operators over $K$ by $K[\Delta]$ and set
$${[\Delta]}:=\{\delta_1^{k_1}\cdots\delta_m^{k_m}~|~k_1,\ldots,k_m\in\mathbb N\}.$$ Similarly, if $\Delta$ is a family of mutually commuting automorphisms of $K$, we set $${[\Delta^{\ast}]}:=\{\delta_1^{k_1}\cdots\delta_m^{k_m}~|~k_1,\ldots,k_m\in\mathbb Z\}$$ and denote the ring of $\Delta^{\ast}$-operators over $K$ by $K[\Delta^{\ast}]$. Furthermore, a free left $K[\Delta]$- (respectively, $K[\Delta^{\ast}]$-) module with a set of free generators $E = \{e_{1},\dots, e_{q}\}$ will be denoted by $K[\Delta]E$ (respectively, $K[\Delta^{\ast}]E$) and  $[\Delta]E$ (respectively, $[\Delta^{\ast}]E$) will denote the set of all elements of the form $\lambda e_{i}$ where $1\leq i\leq q$ and $\lambda\in [\Delta]$ (respectively, $\lambda\in [\Delta^{\ast}]$. Such elements of the free module are called {\em terms}.

Let $E=\{e_1,\ldots,e_q\}$ be a finite set of free generators of the left $K[\Delta^\ast]$ module $K[\Delta^\ast]E$ and $F\subseteq K[\Delta^\ast]E$ finite. There are two popular approaches for computing a Gr\"obner basis of the left $K[\Delta^\ast]$-module $_{K[\Delta^\ast]}\langle F\rangle$. The first approach is due to the second author \cite{Levin5,Levin7} -- the idea also appears \cite{LW} -- works by introducing new variables for the inverses $\delta_1^{-1},\ldots,\delta_m^{-1}$ of $\delta_1,\ldots,\delta_m$ and doing computations in the resulting free module of difference operators. The second approach, originated by the second author \cite{Levin2}, was enhanced by Winkler and Zhou \cite{ZW1,ZW2} who introduced the concept of so-called generalized term orders therefore making $K[\Delta^\ast]E$ a well-ordered set. In the following we will outline the first approach. Proofs for termination and correctness of the algorithms can be found, e.g, in \cite{LW,Levin5,Levin7}.

\subsection{Computing Gr\"obner bases of inversive difference modules via standard bases of associated difference modules}\label{s1}

Let $K$ be a differential (respectively, difference) field with basic set $\Delta=\{\delta_1,\ldots,\delta_m\}$ of derivations (respectively,  endomorphisms) of $K$.

\begin{definition}
 Let $\prec$ be a total order on the set of terms $[\Delta]E$ such that for all elements $1\neq\lambda,\eta,\mu\in[\Delta], e,e'\in E$ we have\begin{enumerate}\item $e\prec\lambda e$, \item $\mu\lambda e\prec\mu\eta e'$ whenever $\lambda e\prec \eta e'$.\end{enumerate}Then the relation $\prec$ is called an \emph{admissible order}.
\end{definition}

For any $f=a_1f_1+\ldots+a_nf_n\in K[\Delta]E$ with $a_1,\ldots,a_n\in K,f_1,\ldots,f_n\in[\Delta]E$ and for a given admissible order $\prec$ we denote the highest term appearing in $f$ with nonzero coefficient by $\lt(f)$, i.e.,$$\lt(f):=\max_\prec\{f_i~|~1\leq i\leq n,a_i\neq 0\},$$and call it the {\em leading term} of $f$. The corresponding coefficient is called the {\em leading coefficient} of $f$ and is denoted by $\lc(f)$.

\begin{definition}
 Let $f,g\in K[\Delta]E\setminus\{0\}$ and $\prec$ an admissible order. If there exists $\lambda\in[\Delta]$ with $\lt(\lambda g)=\lt(f)$ we say that $f$ is \emph{reducible} to $h:=f-\lambda\frac{\lc(f)}{\lc(g)}g$ modulo $g$ and write$$f\longrightarrow_gh.$$Otherwise we say that $f$ is irreducible modulo $g$. Let $G=\{g_1,\ldots,g_p\}\subseteq K[\Delta]E\setminus\{0\}$. If there exist $n\in\mathbb N,f_0,\ldots,f_n,i_1,\ldots,i_n\in\{1,\ldots,p\}$ such that$$f=f_0\longrightarrow_{g_{i_1}}f_1\longrightarrow_{g_{i_2}}\cdots\longrightarrow_{g_{i_n}}f_n=:r$$we say that $f$ is reducible to $r$ modulo $G$. Otherwise we say that $f$ is irreducible modulo $G$.
\end{definition}

The process of reduction is described by the following algorithm.

\begin{algorithm}
\caption{\texttt{Reduction\_algorithm}\label{alg1}}
\begin{algorithmic}
\REQUIRE $0\neq f\in K[\Delta]E$, finite $G\subseteq K[\Delta]E$, and an admissible order $\prec$,
\ENSURE $r\in K[\Delta]E$ such that $f$ is reducible to $r$ modulo $G$ and $r$ is irreducible modulo $G$.
\STATE $r:=f$
\WHILE{there exist $g\in G$ and $\lambda\in [\Delta]$ such that $\lt(\lambda g)=\lt(r)$}
\STATE $r:=r-\lambda\frac{\lc(r)}{\lc(g)}g$
\ENDWHILE
\RETURN $r$
\end{algorithmic}
\end{algorithm}

\begin{definition}
 Let $\prec$ be an admissible order, $N$ a submodule of $K[\Sigma]E$ and $G\subseteq N\setminus\{0\}$ finite such that every $0\neq f\in N$ is reducible to $0$ modulo $G$. Then $G$ is called a \emph{Gr\"obner basis} of the module $N$.
\end{definition}

Every finitely generated $K[\Delta]$-module $M$ has a Gr\"obner basis that can be computed, e.g., via Buchberger's algorithm starting with any finite generating set $\tilde G$  of $M$ (see Algorithm 2 below).

\begin{definition}
 Let $\prec$ be an admissible order on $[\Delta]E$, $g_1,g_2\in K[\Delta]E\setminus\{0\}$, $\lambda_1,\lambda_1\in[\Delta],e_1,e_2\in E$ such that $\lt_\prec(g_1)=\lambda_1 e_1$ and $\lt_\prec(g_2)=\lambda_2 e_2$. The \emph{least common multiple} $\lcm(\lt_\prec(g_1),\lt_\prec(g_2))$ of $\lt_\prec(g_1)$ and $\lt_\prec(g_2)$ is defined by$$\lcm(\lt_\prec(g_1),\lt_\prec(g_2)):=
\begin{cases}
 \lcm(\lambda_1,\lambda_2)e_1&\textnormal{ if }e_1=e_2,\\
0&\textnormal{ if }e_1\neq e_2.
\end{cases}$$
 Let $u_1,u_2\in[\Delta]$ be given by$$u_1:=\frac{\lcm(\lt_\prec(g_1),\lt_\prec(g_2))}{\lt_\prec(g_1)}\qquad\textnormal{ and }\qquad u_2:=\frac{\lcm(\lt_\prec(g_1),\lt_\prec(g_2))}{\lt_\prec(g_2)}.$$Then the \emph{S-polynomial} $S(g_1,g_2)$ of $g_1$ and $g_2$ is defined  by$$S(g_1,g_2):=u_1\frac{g_1}{\lc_\prec(g_1)}-u_2\frac{g_2}{\lc_\prec(g_2)}.$$
\end{definition}

\begin{algorithm}\label{alg2}
\caption{\texttt{Buchberger's\_algorithm}}
\begin{algorithmic}
\REQUIRE $\tilde G\subseteq K[\Delta]E\setminus\{0\}$ finite, $\prec$ an admissible order,
\ENSURE $G\subseteq K[\Delta]E\setminus\{0\}$ being a Gr\"obner basis of $_{K[\Delta]}\langle \tilde G\rangle$.
\STATE $G:=\tilde G$
\WHILE{there exists $g,g'\in G$ such that $S(g,g')$ is not reducible to $0$ modulo $G$}
\STATE $G:=G\cup\{\texttt{Reduction\_algorithm}(S(g,g'),G, \prec)\}$
\ENDWHILE
\RETURN $G$
\end{algorithmic}
\end{algorithm}

The following theorem being a special case of \cite[Thm. 4.12]{Levin4} and \cite[Thm. 3.3.15]{Levin7} describes how to obtain the dimension polynomial associated with a system of differential or difference equations.

\begin{theorem}\label{thm1}
 Let $M$ be a difference $K$-vector space generated (as a left $K[\Delta]$-module) by elements $m_1,\ldots,m_q$, $F$ a free $K[\Delta]$-module with set of free generators $E=\{e_1,\ldots,e_q\}$, $\pi:F\to M$ the difference epimorphism $(e_i\mapsto m_i$ for $i=1,\dots, q$) and $N:=\ker(\pi)$. Let $G\subseteq K[\Delta]E$ be a Gr\"obner basis of $N$ with respect to the term order $\prec$ defined by\begin{eqnarray*}\lefteqn{\delta_1^{k_1}\cdots\delta_m^{k_m}e_i\prec\delta_1^{l_1}\cdots\delta_m^{l_m}e_j}\nonumber\\&\Longleftrightarrow&(k_1+\cdots+k_m,i,k_1,\ldots,k_m)<_{\lex}(l_1+\cdots+l_m,j,l_1,\ldots,l_m),\nonumber\end{eqnarray*}where $<_{\lex}$ denotes the lexicographic order. For $r\in\mathbb Z$ let\begin{eqnarray*}M_r&:=&\{\lambda m\in[\Delta]\{m_1,\ldots,m_q\}~|~\ord\lambda\leq r\}\textnormal{ and}\nonumber\\U_r&:=&\{\lambda e\in [\Delta]E~|~\ord\lambda\leq r\,\,\,  \textnormal{and}\,\,\, \lambda e\neq \mu\lt_\prec(g)\,\,\, \textnormal{for any}\,\,\, g\in G, \mu\in [\Delta] \}.\nonumber\end{eqnarray*}Then $(M_r)_{r\in\mathbb Z}$ is an excellent filtration of $M$ and for any $r\in \mathbb N$ the set $\pi(U_r)$ is a basis for the $K$-vector space $M_r$
\end{theorem}
The following proposition is obtained from \cite[Prop. 2.2.11.]{KLMP} by realizing that a term $\lambda e$ is irreducible if and only if there exist no $\eta\in[\Delta],g\in G$ with $\lt(g)=\mu e$ and $\eta\mu=\lambda$.

\begin{proposition}\label{prop1}
 With the notation of Theorem \ref{thm1} for every $i=1,\dots, q$, let $$G_i:=\{\lt(g)~|~g\in G,~\lt(g)\in[\Delta]e_i\}$$ and let $\Lambda_i=(\lambda_{i,j,k})\in\mathbb N^{|G_i|\times 2m}$ satisfy the following condition:  for every$$\lt(g)=\alpha_1^{a_1}\cdots\alpha_m^{a_m}\beta_1^{b_1}\cdots\beta_m^{b_m}e_i\in G_i$$there exists $j\in\{1,\ldots,|G_i|\}$ with $(\lambda_{i,j,1},\ldots,\lambda_{i,j,2m})=(a_1,\ldots,a_m,b_1,\ldots,b_m)$.

 Furthermore, for any $l,n\in\mathbb N$ with $1\leq n$ and $0\leq l\leq n$, let $A(l, n)$ denote the set of all $l$-element subsets of $\{1,\ldots,n\}$ and for every $1\leq i\leq q,\emptyset\neq\xi\in A(l,|G_i|)$, let $\lambda_{i,\xi,k}:=\max_{j\in\xi} \lambda_{i,j,k}$ and $\lambda_{i,\emptyset,k}:=0$. Finally, for any $1\leq i\leq q,\xi\in A(l,|G_i|)$ let $f_{i,\xi}:=\sum_{k=1}^{2m}\lambda_{i,\xi,k}$. Then$$|U_r|=\sum_{i=1}^r\sum_{l=0}^{|G_i|}(-1)^l\sum_{\xi\in A(l,|G_i|)}\binom{r+2m-f_{i,\xi}}{2m}.$$
\end{proposition}

An idea of the following kind was also by Ziming Li and Min Wu \cite{LW}. Let $K$ be an inversive difference field with basic set of automorphisms $\Delta=\{\delta_1,\ldots,\delta_m\}$ and $F\subseteq K[\Delta^\ast]E$ a finite set of generators for a left $K[\Delta^\ast]$ module. Let $\alpha_1,\ldots,\alpha_m,\beta_1,\ldots,\beta_m$ be endomorphisms of $K$ such that for $i=1,\ldots,m$ and $k\in K$ we have$$\alpha_i(k):=\delta_i(k),\qquad\textnormal{ and }\qquad\beta_i(k):=\delta_i^{-1}(k).$$Then $K$ can be considered as a difference field with basic set $\Sigma:=\{\alpha_1,\ldots,\alpha_m,$ $\beta_1,\ldots,$ $\beta_m\}$. By $\rho:K[\Delta^\ast]\to K[\Sigma]/_{K[\Sigma]}\langle\{\alpha_i\beta_i e-e~|~1\leq i\leq m,e\in E\}\rangle$ we denote the natural isomorphism 
$$\rho:\delta_{1}^{k_1}\dots\delta_{m}^{k_m}e\mapsto\alpha_1^{\max\{k_1,0\}}\dots\alpha_{m}^{\max\{k_m,0\}}\beta_1^{\max\{-k_1,0\}}\dots\beta_{m}^{\max\{-k_m,0\}}e\qquad(e\in E).$$
Let $\tilde F:=\rho(F)\cup\{\alpha_i\beta_ie-e~|~1\leq i\leq m,~e\in E\}$. Then, $K[\Delta^\ast]E/_{K[\Delta^\ast]}\langle F\rangle$ is isomorphic to $K[\Sigma]E/_{K[\Sigma]}\langle \tilde F\rangle$, so in order to compute the $\Delta^{\ast}$-dimension polynomial of a $\Delta^{\ast}$-$K$-module $K[\Delta^\ast]E/_{K[\Delta^\ast]}\langle F\rangle$ associated with a finite system of generators $F$, it suffices to compute the $\Sigma$-dimension polynomial of the $\Sigma$-$K$module $K[\Sigma]E/_{K[\Sigma]}\langle \tilde F\rangle$ associated with the set of generators $\tilde F$.

\subsection{Examples for the computation of differential and difference dimension polynomials}

In this subsection we give several examples for the computation of differential dimension polynomials associated with systems of differential equations arising in mathematical physics and of difference dimension polynomials associated with their difference schemes.

\begin{example}\textnormal{\textbf{(Diffusion equation in $1$-space)}}\\
 The diffusion equation in one spatial dimension for a constant collective diffusion coefficient $a$ and unknown function $ u(x,t)$ describing the density of the diffusing material at given position $x$ and time $t$ is given by\begin{equation}\label{diffusionequation}\frac{\partial u(x,t)}{\partial t}=a\frac{\partial^2 u(x,t)}{\partial x^2}.\end{equation}

 \centerline{\textnormal{\bf{Differential dimension polynomial}}}Let $K$ be a differential field with basic set $\Delta=\{\delta_x=\frac{\partial}{\partial x},\delta_t=\frac{\partial}{\partial t}\}$ containing $a$ and let $M$ be a differential $K$-vector space generated as $K[\Delta]$-module by one generator $m$ satisfying the defining equation$$\delta_tm=a\delta_x^2m.$$Then $M$ is isomorphic to the factor module of a free $K[\Delta]$-module with free generator $e$ by its submodule $N$ generated by$$G:=\{\delta_te-a\delta_x^2e\}.$$Since $G$ consists of only one element there are no S-polynomials. Therefore $G$ is already a Gr\"obner basis of $N$ for any admissible order on $[\Delta]e$. Let the admissible order $\prec$ on $[\Delta]e$ be given by$$\delta_x^{k_x}\delta_t^{k_t}e\prec\delta_x^{l_x}\delta_t^{l_t}e:\Longleftrightarrow(k_x+k_t,k_x,k_t)<_{\lex}(l_x+l_t,l_x,l_t).$$Then for all $2\leq r\in\mathbb N$ we have\begin{eqnarray*}U_r&=&\{\delta_x^{k_x}\delta_t^{k_t}e~|~k_x+k_t\leq r,~\delta_x^{k_x}\delta_t^{k_t}e\textnormal{ is irreducible modulo }G\}\nonumber\\&=&\{e,\delta_te,\ldots,\delta_t^re,\delta_xe,\delta_x\delta_te,\ldots,\delta_x\delta_t^{r-1}e\},\nonumber\end{eqnarray*}and therefore $|U_r|=2r+1$.

 Thus, the differential dimension polynomial associated with the diffusion equation in one spatial dimension for a constant collective diffusion coefficient is given by $\phi(r)=2r+1$.

\centerline{\textnormal{\bf{Difference dimension polynomial for forward difference scheme}}}
In order to obtain a forward difference scheme for the diffusion equation \eqref{diffusionequation} every occurence of $\frac{\partial u(x,t)}{\partial x}$ and $\frac{\partial u(x,t)}{\partial t}$ is replaced by $u(x+1,t)- u(x,t)$ and $ u(x,t+1)- u(x,t)$, respectively. We obtain
\begin{equation}\label{forwarddifferencediffusionequation} u(x,t+1)- u(x,t)=a(u(x+2,t)-2u(x+1,t)+ u(x,t)).\end{equation}
Let $K$ be an inversive difference field with basic set $\Delta=\{\delta_x:x\mapsto x+1,\delta_t:t\mapsto t+1\}$ containing $a$ and let $M$ be an inversive difference $K$-vector space generated as a left $K[\Delta^\ast]$-module by one generator $m$ satisfying the defining equation
$$\delta_tm-m=a(\delta_x^2m-2\delta_xm+m).$$
Then $M$ is isomorphic to the factor module of a free $K[\Delta^\ast]$-module with free generator $e$ by its submodule $N$ generated by$$G:=\{\delta_te-a\delta_x^2e+2a\delta_xe-(1+a)e\}.$$We will compute the difference dimension polynomial associated with the difference scheme \eqref{forwarddifferencediffusionequation} using the method described at the end of subsection A. Thus, we consider $K$ as a difference field with basic set $\Sigma=\{\alpha_x:x\mapsto x+1,\alpha_t:t\mapsto t+1,\beta_x:x\mapsto x-1,\beta_t:t\mapsto t-1\}$. Let $\tilde G:=\{g_1:=\alpha_te-a\alpha_x^2e+2a\alpha_xe-(1+a)e,g_2:=\alpha_x\beta_xe-e,g_3:=\alpha_t\beta_te-e\}$ and $I=_{K[\Sigma]}\langle \tilde G\rangle$. Then $K[\Sigma]e/I$ is isomorphic $K[\Delta^\ast]e/N$ via the isomorphism
$$\alpha_x^{a_x}\alpha_t^{a_t}\beta_x^{b_x}\beta_t^{b_t}e\mapsto\delta_x^{a_x-b_x}\delta_t^{a_t-b_t}e.$$We fix an admissible order $\prec$ on $[\Sigma]e$ defined by\begin{eqnarray*}\lefteqn{\alpha_x^{a_x}\alpha_t^{a_t}\beta_x^{b_x}\beta_t^{b_t}e\prec\alpha_x^{c_x}\alpha_t^{c_t}\beta_x^{d_x}\beta_t^{d_t}e:\Longleftrightarrow}\nonumber\\&&(a_x+a_t+b_x+b_t,a_x,a_t,b_x,b_t)<_{\lex}(c_x+c_t+d_x+d_t,c_x,c_t,d_x,d_t)\nonumber\end{eqnarray*}
and compute a Gr\"obner basis of $I$ with respect to $\prec$. The S-polynomial $S(g_1,g_2)$ of $g_1$ and $g_2$ is given by$$S(g_1,g_2)=-2 \alpha_x \beta_xe- \frac{1}{a}\alpha_t\beta_xe+\left(1+\frac{1}{a}\right) \beta_xe+\alpha_xe$$ and is reducible modulo $g_2$ to $ - \frac{1}{a}\alpha_t\beta_xe+\left(1+\frac{1}{a}\right) \beta_xe+\alpha_xe-2e$, which is irreducible modulo $\tilde G$. Hence, $g_4:=- \frac{1}{a}\alpha_t\beta_xe+\left(1+\frac{1}{a}\right) \beta_xe+\alpha_xe-2e$ should be inserted into $\tilde G$. For $S(g_1,g_3)=   -2 \alpha_x\alpha_t \beta_te - \frac{1}{a}\alpha_t^2\beta_te+\left(1+\frac{1}{a}\right) \alpha_t \beta_te+\alpha_x^2e$ we have\begin{eqnarray*}S(g_1,g_3)&\longrightarrow_{g_3}& -\frac{1}{a}\alpha_t^2 \beta_t e+\alpha_x^2 e+\left(1+\frac{1}{a}\right) \alpha_t \beta_t e-2 \alpha_x e\nonumber\\&\longrightarrow_{g_3}&\alpha_x^2 e+\left(1+\frac{1}{a}\right) \alpha_t \beta_t e-2 \alpha_x e-\frac{1}{a}\alpha_t e\nonumber\\&\longrightarrow_{g_1}&\left(1+\frac{1}{a}\right) \alpha_t \beta_t e-\left(1+\frac{1}{a}\right) e\nonumber\\&\longrightarrow_{g_3}&0.\nonumber\end{eqnarray*}Furthermore\begin{eqnarray*}
S(g_1,g_4)&=&a \alpha_x^3 e+(a+1) \alpha_x^2 \beta_x e-2 \alpha_x \alpha_t \beta_x e-\frac{1}{a}\alpha_t^2 \beta_x e-2 a \alpha_x^2 e\nonumber\\&&+\left(1+\frac{1}{a}\right) \alpha_t \beta_x e\nonumber\\&\longrightarrow_{g_1}& (a+1) \alpha_x^2 \beta_x e-2 \alpha_x \alpha_t \beta_x e-\frac{1}{a}\alpha_t^2 \beta_x e+\alpha_x \alpha_t e\nonumber\\&&+\left(1+\frac{1}{a}\right) \alpha_t \beta_x e-(1+a) \alpha_x e\nonumber\\&\longrightarrow_{g_1}& -2 \alpha_x \alpha_t \beta_x e-\frac{1}{a}\alpha_t^2 \beta_x e+\alpha_x \alpha_t e+(2 a+2) \alpha_x \beta_x e\nonumber\\&&+(2+2/a) \alpha_t \beta_x e-(1+a) \alpha_x e+(-2-a-1/a) \beta_x e\nonumber\\&
\longrightarrow_{g_2}& -\frac{1}{a}\alpha_t^2 \beta_x e+\alpha_x \alpha_t e+(2 a+2) \alpha_x \beta_x e+(2+2/a) \alpha_t \beta_x e\nonumber\\&&-(1+a) \alpha_x e-2 \alpha_t e+(-2-a-1/a) \beta_x e\nonumber\\&\longrightarrow_{g_4}& (2 a+2) \alpha_x \beta_x e+\left(1+\frac{1}{a}\right) \alpha_t \beta_x e-(1+a) \alpha_x e\nonumber\\&&+(-2-a-1/a) \beta_x e\nonumber\\&\longrightarrow_{g_2}& \left(1+\frac{1}{a}\right) \alpha_t \beta_x e-(1+a) \alpha_x e+(-2-a-1/a) \beta_x e\nonumber\\&&+(2 a+2) e\nonumber\\&\longrightarrow_{g_4}& 0\nonumber\\
S(g_2, g_3)&=& \alpha_x \beta_x e-\alpha_t \beta_t e\nonumber\\&\longrightarrow_{g_2}& -\alpha_t \beta_t e+e\nonumber\\&\longrightarrow_{g_3}& 0\nonumber\\
S(g_2,g_4)&=& a \alpha_x^2 e+(a+1) \alpha_x \beta_x e-2 a \alpha_x e-\alpha_t e\nonumber\\&\longrightarrow_{g_1}& (a+1) \alpha_x \beta_x e-(1+a) e\nonumber\\&\longrightarrow_{g_2}& 0\nonumber\\
S(g_3,g_4)&=& a \alpha_x \beta_t e+(a+1) \beta_x \beta_t e-\beta_x e-2 a \beta_t e\nonumber\\&=:&g_5\nonumber\\
S(g_1, g_5)&=& -\left(1+\frac{1}{a}\right) \alpha_x \beta_x \beta_t e+\frac{1}{a}\alpha_x \beta_x e-\frac{1}{a}\alpha_t \beta_t e\nonumber\\&&+\left(1+\frac{1}{a}\right) \beta_t e\nonumber\\&\longrightarrow_{g_2}& \frac{1}{a}\alpha_x \beta_x e-\frac{1}{a}\alpha_t \beta_t e\nonumber\\&\longrightarrow_{g_2}& -\frac{1}{a}\alpha_t \beta_t e+\frac{1}{a}e\nonumber\\&\longrightarrow_{g_3}& 0\nonumber\\
S(g_2, g_5)&=& -\left(1+\frac{1}{a}\right) \beta_x^2 \beta_t e+\frac{1}{a}\beta_x^2 e+2 \beta_x \beta_t e-\beta_t e\nonumber\\&=:&g_6\nonumber
\end{eqnarray*}

Further computations of the $S$-polynomials $S(g_{i}, g_{6})$, $1\leq i\leq 5$, show that all of them are reducible to $0$ modulo $\{g_{1},\dots, g_{6}\}$. Hence, a Gr\"obner basis for $I$ is given by\begin{eqnarray*}\Biggl\{
g_1&=&\alpha_te-a\alpha_x^2e+2a\alpha_xe-(1+a)e,
\nonumber\\g_2&=&\alpha_x\beta_xe-e,
\nonumber\\g_3&=&\alpha_t\beta_te-e,
\nonumber\\g_4&=&-\frac{1}{a}\alpha_t\beta_xe+\left(1+\frac{1}{a}\right) \beta_xe+\alpha_xe-2e,
\nonumber\\g_5&=&a \alpha_x \beta_t e+(a+1) \beta_x \beta_t e-\beta_x e-2 a \beta_t e,
\nonumber\\g_6&=&-\left(1+\frac{1}{a}\right) \beta_x^2 \beta_t e+\frac{1}{a}\beta_x^2 e+2 \beta_x \beta_t e-\beta_t e\Biggr\}.\nonumber\end{eqnarray*}Applying Proposition \ref{prop1} we obtain $$|\{\lambda\in[\Sigma]e~|~\ord(\lambda)\leq r,\lambda\textnormal{ is irreducible modulo }G\}|=5r.$$ for all sufficiently large  $r$. Therefore, the inversive difference dimension polynomial associated with the difference scheme \eqref{forwarddifferencediffusionequation} is given by $\phi(r)=5r$.

\centerline{\textnormal{\bf Difference dimension polynomial for symmetric difference scheme}} In order to obtain a space symmetric difference scheme for the diffusion equation \eqref{diffusionequation} every occurrence of $\D\frac{\partial^2 u(x,t)}{\partial x^2}$ and\, $\D\frac{\partial u(x,t)}{\partial t}$ is replaced by $u(x+1,t)-2 u(x,t)+u(x-1,t)$ and $ u(x,t+1)- u(x,t)$, respectively. We obtain\begin{equation}\label{symmetricdifferencediffusionequation} u(x,t+1)- u(x,t)=a(u(x+1,t)-2 u(x,t)+u(x-1,t)).\end{equation}Let $K$ be an inversive difference field with basic set $\Delta=\{\delta_x:x\mapsto x+1,\delta_t:t\mapsto t+1\}$ containing $a$ and let $M$ be an inversive difference $K$-vector space generated as a left $K[\Delta^\ast]$-module by one generator $m$ satisfying the defining equation$$\delta_tm-m=a(\delta_xm-2m+\delta_x^{-1}m).$$Then $M$ is isomorphic to the factor module of a free $K[\Delta^\ast]$-module with free generator $e$ by its submodule $N$ generated by$$G:=\{\delta_te-a\delta_xe-\delta_x^{-1}e+(2a-1)e\}.$$Now consider $K$ as a difference field with basic set $\Sigma=\{\alpha_x:x\mapsto x+1,\alpha_t:t\mapsto t+1,\beta_x:x\mapsto x-1,\beta_t:t\mapsto t-1\}$. Let $\tilde G:=\{g_1:=\alpha_te-\alpha_xe-\beta_xe+(2a-1)e,g_2:=\alpha_x\beta_xe-e,g_3:=\alpha_t\beta_te-e\}$ and $I=_{K[\Sigma]}\langle \tilde G\rangle$. Then $K[\Sigma]e/I$ is isomorphic $K[\Delta^\ast]e/N$ via the isomorphism$$\alpha_x^{a_x}\alpha_t^{a_t}\beta_x^{b_x}\beta_t^{b_t}e\mapsto\delta_x^{a_x-b_x}\delta_t^{a_t-b_t}e.$$ Let us fix an admissible order $\prec$ on $[\Sigma]e$ defined by\begin{eqnarray*}\lefteqn{\alpha_x^{a_x}\alpha_t^{a_t}\beta_x^{b_x}\beta_t^{b_t}e\prec\alpha_x^{c_x}\alpha_t^{c_t}\beta_x^{d_x}\beta_t^{d_t}e:\Longleftrightarrow}\nonumber\\&&(a_x+a_t+b_x+b_t,a_x,a_t,b_x,b_t)<_{\lex}(c_x+c_t+d_x+d_t,c_x,c_t,d_x,d_t)\nonumber\end{eqnarray*}A Gr\"obner basis of $I$ is then given by
\begin{eqnarray*}
\{g_1&:=& a \beta_x^2 \beta_te-(2 a-1) \beta_x \beta_te-\beta_x+a \beta_te,\nonumber\\g_2&:=& -\frac{1}{a}\alpha_t \beta_xe+\beta_x^2e-\left(2-\frac{1}{a}\right) \beta_xe+e,\nonumber\\g_3&:=& \alpha_t \beta_te-e,\nonumber\\g_4&:=& a \alpha_xe-\alpha_te+a \beta_xe-(2 a-1)e\}.\nonumber\end{eqnarray*}Applying Proposition \ref{prop1} we obtain $$|\{\lambda\in[\Sigma]e~|~\ord(\lambda)\leq r,\lambda\textnormal{ is irreducible modulo }G\}|=4r$$ for all sufficiently large $r$. Hence, the inversive difference dimension polynomial associated with the difference scheme \eqref{symmetricdifferencediffusionequation} is given by $\phi(r)=4r$.

Thus, the symmetric difference scheme for the diffusion equation has higher strength (that is, smaller dimension polynomial) than the forward scheme, so the symmetric scheme is more preferable from the point of view of strength.

\end{example}

\begin{example}\textnormal{\textbf{(Maxwell equations for vanishing free current density and  free charge density)}}\\
 Let $E=(E_1,E_2,E_3),\, D=(D_1,D_2,D_3),\, H=(H_1,H_2,H_3),\, B=(B_1,B_2,B_3),$ $J_f=(J_1,J_2,J_3)$ and $\rho_f$ be functions of $(x,y,z,t)$ that denote  electric field strength, electric displacement vector, magnetic field strength, magnetic displacement vector, free current density and free charge density, respectively. With$$\nabla:=\left(\frac{\partial}{\partial x},\frac{\partial}{\partial y},\frac{\partial}{\partial z}\right)$$
Maxwell's equations in $3$ spatial dimensions are given by$$\nabla\cdot D=\rho_f,\qquad\nabla\cdot B=0,\qquad\nabla\times E+\frac{\partial B}{\partial t}=0,\quad\textnormal{and}\quad \nabla\times H=J_f+\frac{\partial D}{\partial t}.$$ Assuming $J_f=0$ and $\rho_f=0$, Maxwell's equations can be considered as a set of homogeneous linear differential equations.

\centerline{\textnormal{\bf{Differential dimension polynomial}}}
Let $K$ be a differential field with basic set $\Delta = \{\delta_x=\frac{\partial}{\partial x},\delta_y=\frac{\partial}{\partial y},\delta_z=\frac{\partial}{\partial z},\delta_t=\frac{\partial}{\partial t}\}$. Assuming $J_f=0$ and $\rho_f=0$ Maxwell's equations give rise to a differential $K[\Delta]$-module $M$ with generators $e_1,e_2,e_3,d_1,d_2,d_3,h_1,h_2,h_3,b_1,b_2,b_3$ satisfying
\begin{eqnarray*}\delta_xd_1+\delta_yd_2+\delta_zd_3&=~~0~~=&\delta_xb_1+\delta_yb_2+\delta_zb_3,\nonumber\\\delta_ye_3-\delta_ze_2+\delta_tb_1&=~~0~~=&\delta_yh_3-\delta_zh_2-\delta_td_1,\nonumber\\\delta_ze_1-\delta_xe_3+\delta_tb_2&=~~0~~=&\delta_zh_1-\delta_xh_3-\delta_td_2,\nonumber\\\delta_xe_2-\delta_ye_1+\delta_tb_3&=~~0~~=&\delta_xh_2-\delta_yh_1-\delta_td_3.\nonumber\end{eqnarray*}Then $M$ is isomorphic to the factor module of a free $K[\delta_x,\delta_y,\delta_z,\delta_t]$-module with free generators $p_1,\ldots,p_{12}$ by its submodule $N$ generated by
\begin{eqnarray*}G&=&\{\delta_xp_7+\delta_yp_8+\delta_zp_9,\delta_xp_{10}+\delta_yp_{11}+\delta_zp_{12},\delta_yp_3-\delta_zp_2+\delta_tp_{10},\nonumber\\&&~\delta_yp_{6}-\delta_zp_{5}-\delta_tp_7,\delta_zp_1-\delta_xp_3+\delta_tp_{11},\delta_zp_{4}-\delta_xp_{6}-\delta_tp_8,\nonumber\\&&~\delta_xp_2-\delta_yp_1+\delta_tp_{12},\delta_xp_{5}-\delta_yp_{4}-\delta_tp_9\}.\nonumber\end{eqnarray*}
We define an admissible order $\prec$ by
\begin{eqnarray*}\lefteqn{\delta_x^{a_x}\delta_y^{a_y}\delta_z^{a_z}\delta_t^{a_t}e_{j_1}\prec\delta_x^{b_x}\delta_y^{b_y}\delta_z^{b_z}\delta_t^{b_t}e_{j_2}:\Longleftrightarrow}\nonumber\\&&(a_x+a_y+a_z+a_t,j_1,a_x,a_y,a_z,a_t)\nonumber\\&&<_{\lex}(b_x+b_y+b_z+b_t,j_2,b_x,b_y,b_z,b_t).\nonumber\end{eqnarray*}
Then $G$ is a Gr\"obner basis and by Proposition \ref{prop1} the differential dimension polynomial associated with Maxwell's equations for vanishing free current density and free charge density is given by$$\phi(r)=\frac{1}{4}r^4+\frac{19}{6}r^3+\frac{55}{4}r^2+\frac{137}{6}r+12.$$

\centerline{\textnormal{\bf{Difference dimension polynomial for forward difference scheme}}}
Let $K$ be an inversive difference field with basic set $\Delta=\{\delta_x:x\mapsto x+1,\delta_y:y\mapsto y+1,\delta_z:z\mapsto z+1,\delta_t:t\mapsto t+1\}$. If we replace every occurrence of $\delta_{i}$ in $G$ by $\delta_{i}-1$, where $i\in \{x, y, z, t\}$, we obtain a set
\begin{eqnarray*}\tilde G&=&\{\delta_{x}p_{7}-p_{7}+\delta_{y}p_{8}-p_{8}+\delta_{z}p_{9}-p_{9},\nonumber\\&&~\delta_{x}p_{1}-p_{10}+\delta_{y}p_{1}-p_{11}+\delta_{z}p_{2}-p_{12},\nonumber\\&&~\delta_{y}p_{3}-p_{3}-\delta_{z}p_{2}-p_{2}+\delta_{t}p_{1}-p_{10},\nonumber\\&&~\delta_{y}p_{6}-p_{6}-\delta_{z}p_{5}-p_{5}-\delta_{t}p_{7}-p_{7},\nonumber\\&&~\delta_{z}p_{1}-p_{1}-\delta_{x}p_{3}-p_{3}+\delta_{t}p_{1}-p_{11},\nonumber\\&&~\delta_{z}p_{4}-p_{4}-\delta_{x}p_{6}-p_{6}-\delta_{t}p_{8}-p_{8},\nonumber\\&&~\delta_{x}p_{2}-p_{2}-\delta_{y}p_{1}-p_{1}+\delta_{t}p_{2}-p_{12},\nonumber\\&&~\delta_{x}p_{5}-p_{5}-\delta_{y}p_{4}-p_{4}-\delta_{t}p_{9}-p_{9}\}.\nonumber\end{eqnarray*}
associated with Maxwell's equations for vanishing free current density and free charge density. Applying the Buchberger algorithm we obtain the following 80-element Gr\"obner basis $G$ of the associated $K[\alpha_x,\alpha_y,\alpha_z,$ $\alpha_t,\beta_x,\beta_y,\beta_z,\beta_t]$-submodule of the free module with free generators $p_1,\ldots,$ $p_{12}$.
\begin{eqnarray*}
 &G=\{ \beta_x \beta_y \beta_z p_{12}+\beta_x \beta_y \beta_z p_{11}-\beta_x \beta_y p_{12}+\beta_x \beta_y \beta_z p_{10}-\beta_x \beta_z p_{11}-\beta_y \beta_z p_{10},
 \beta_y \beta_z p_{12}-\alpha_x \beta_y \beta_z p_{10}&\\&
+\beta_y \beta_z p_{11}-\beta_y p_{12}+\beta_y \beta_z p_{10}-\beta_z p_{11},
 \beta_x \beta_y \beta_z p_{9}+\beta_x \beta_y \beta_z p_{8}-\beta_x \beta_y p_{9}+\beta_x \beta_y \beta_z p_{7}&\\&
-\beta_x \beta_z p_{8}-\beta_y \beta_z p_{7},
 \beta_y \beta_z p_{9}-\alpha_x \beta_y \beta_z p_{7}+\beta_y \beta_z p_{8}-\beta_y p_{9}+\beta_y \beta_z p_{7}-\beta_z p_{8},
 \beta_y \beta_z \beta_t p_{7}&\\&
-\beta_y \beta_z \beta_t p_{6}-\beta_y \beta_z p_{7}+\beta_y \beta_z \beta_t p_{5}+\beta_z \beta_t p_{6}-\beta_y \beta_t p_{5},
 -\beta_y \beta_t p_{7}+\alpha_z \beta_y \beta_t p_{5}+\beta_y \beta_t p_{6}+\beta_y p_{7}&\\&
-\beta_y \beta_t p_{5}-\beta_t p_{6},
 -\beta_x \beta_z \beta_t p_{8}+\beta_x \beta_z p_{8}-\beta_x \beta_z \beta_t p_{6}+\beta_z \beta_t p_{6}+\beta_x \beta_z \beta_t p_{4}-\beta_x \beta_t p_{4},
 \beta_x \beta_t p_{8}&\\&
-\beta_x p_{8}+\beta_x \beta_t p_{6}+\alpha_z \beta_x \beta_t p_{4}-\beta_t p_{6}-\beta_x \beta_t p_{4},
 \beta_x \beta_y \beta_t p_{9}-\beta_x \beta_y p_{9}-\beta_x \beta_y \beta_t p_{5}+\beta_x \beta_y \beta_t p_{4}&\\&
+\beta_y \beta_t p_{5}-\beta_x \beta_t p_{4},
 \beta_x \beta_t p_{9}-\beta_x p_{9}-\alpha_y \beta_x \beta_t p_{4}-\beta_x \beta_t p_{5}+\beta_x \beta_t p_{4}+\beta_t p_{5},
 -\beta_y \beta_z \beta_t p_{10}&\\&
+\beta_y \beta_z p_{10}-\beta_y \beta_z \beta_t p_{3}+\beta_y \beta_z \beta_t p_{2}+\beta_z \beta_t p_{3}-\beta_y \beta_t p_{2},
 \beta_y \beta_t p_{10}-\beta_y p_{10}+\alpha_z \beta_y \beta_t p_{2}+\beta_y \beta_t p_{3}&\\&
-\beta_y \beta_t p_{2}-\beta_t p_{3},
 \beta_x \beta_z \beta_t p_{11}-\beta_x \beta_z p_{11}-\beta_x \beta_z \beta_t p_{3}+\beta_z \beta_t p_{3}+\beta_x \beta_z \beta_t p_{1}-\beta_x \beta_t p_{1},
 -\beta_x \beta_t p_{11}&\\&
+\beta_x p_{11}+\beta_x \beta_t p_{3}+\alpha_z \beta_x \beta_t p_{1}-\beta_t p_{3}-\beta_x \beta_t p_{1},
 -\beta_x \beta_y \beta_t p_{12}+\beta_x \beta_y p_{12}-\beta_x \beta_y \beta_t p_{2}&\\&
+\beta_x \beta_y \beta_t p_{1}+\beta_y \beta_t p_{2}-\beta_x \beta_t p_{1},
 -\beta_x \beta_t p_{12}+\beta_x p_{12}-\alpha_y \beta_x \beta_t p_{1}-\beta_x \beta_t p_{2}+\beta_x \beta_t p_{1}+\beta_t p_{2},&\\&
 -\beta_t p_{9}+p_{9}-\alpha_x \beta_t p_{5}+\alpha_y \beta_t p_{4}+\beta_t p_{5}-\beta_t p_{4},
 -\beta_t p_{8}+\alpha_x \beta_t p_{6}+p_{8}-\beta_t p_{6}-\alpha_z \beta_t p_{4}&\\&
+\beta_t p_{4},
 -\alpha_y \beta_t p_{6}-\beta_t p_{7}+\alpha_z \beta_t p_{5}+\beta_t p_{6}+p_{7}-\beta_t p_{5},
 -\beta_t p_{12}+p_{12}+\alpha_x \beta_t p_{2}-\alpha_y \beta_t p_{1}&\\&
-\beta_t p_{2}+\beta_t p_{1},
 -\beta_t p_{11}+p_{11}-\alpha_x \beta_t p_{3}+\beta_t p_{3}+\alpha_z \beta_t p_{1}-\beta_t p_{1},
 -\beta_t p_{10}+p_{10}+\alpha_y \beta_t p_{3}&\\&
-\alpha_z \beta_t p_{2}-\beta_t p_{3}+\beta_t p_{2},
 \alpha_y \beta_z p_{11}-\beta_z p_{12}+\alpha_x \beta_z p_{10}-\beta_z p_{11}+p_{12}-\beta_z p_{10},
 \alpha_y \beta_z p_{8}-\beta_z p_{9}&\\&
+\alpha_x \beta_z p_{7}-\beta_z p_{8}+p_{9}-\beta_z p_{7},
 -\alpha_t p_{9}+p_{9}+\alpha_x p_{5}-\alpha_y p_{4}-p_{5}+p_{4},
 -\alpha_t p_{8}+p_{8}-\alpha_x p_{6}&\\&
+p_{6}+\alpha_z p_{4}-p_{4},
 -\alpha_t p_{7}+\alpha_y p_{6}+p_{7}-\alpha_z p_{5}-p_{6}+p_{5},
 \alpha_t p_{12}-p_{12}+\alpha_x p_{2}-\alpha_y p_{1}-p_{2}&\\&
+p_{1},
 \alpha_t p_{11}-p_{11}-\alpha_x p_{3}+p_{3}+\alpha_z p_{1}-p_{1},
 \alpha_t p_{10}-p_{10}+\alpha_y p_{3}-\alpha_z p_{2}-p_{3}+p_{2},
 \alpha_z p_{12}&\\&
+\alpha_y p_{11}-p_{12}+\alpha_x p_{10}-p_{11}-p_{10},  \alpha_z p_{9}+\alpha_y p_{8}-p_{9}+\alpha_x p_{7}-p_{8}-p_{7}\}&\nonumber\\
&\cup\left\{\alpha_i\beta_ip_j-p_j~\left|~i\in\{x,y,z,t\},j\in\{1,\ldots,12\}\right.\right\}.&\nonumber
\end{eqnarray*}

Applying Proposition \ref{prop1} we obtain that the difference dimension polynomial associated with the forward difference scheme is given by$$\phi(r)=4r^4+18r^3+35r^2+31r+12.$$

\centerline{\textnormal{\bf{Difference dimension polynomial for symmetric difference scheme}}}
~\\Let $K$ be an inversive difference field with basic set $\Delta=\{\delta_x:x\mapsto x+1,\delta_y:y\mapsto y+1,\delta_z:z\mapsto z+1,\delta_t:t\mapsto t+1\}$. Using the symmetric difference scheme we replace every occurrence of $\delta_i$ in $G$  by $\frac{1}{2}(\delta_i-\delta_i^{-1})$ ($i\in\{x,y,z,t\}$) and arrive at a set
\begin{eqnarray*}\tilde
G&=&\{\delta_{x}p_{7}-p_{7}+\delta_{y}p_{8}-p_{8}+\delta_{z}p_{9}-p_{9},\nonumber\\&&~\delta_{x}p_{10}-p_{10}+
\delta_{y}p_{11}-p_{11}+\delta_{z}p_{12}-p_{12},\nonumber\\&&~\delta_{y}p_{3}-p_{3}-\delta_{z}p_{2}-p_{2}+
\delta_{t}p_{10}-p_{10},\nonumber\\&&~\delta_{y}p_{6}-p_{6}-\delta_{z}p_{5}-p_{5}-\delta_{t}p_{7}-p_{7},
\nonumber\\&&~\delta_{z}p_{1}-p_{1}-\delta_{x}p_{3}-p_{3}+\delta_{t}p_{11}-p_{11},
\nonumber\\&&~\delta_{z}p_{4}-p_{4}-\delta_{x}p_{6}-p_{6}-\delta_{t}p_{8}-p_{8},
\nonumber\\&&~\delta_{x}p_{2}-p_{2}-\delta_{y}p_{1}-p_{1}+\delta_{t}p_{12}-p_{12},
\nonumber\\&&~\delta_{x}p_{5}-p_{5}-\delta_{y}p_{4}-p_{4}-\delta_{t}p_{9}-p_{9}\}.
\nonumber\end{eqnarray*}Proceeding as above we obtain that the corresponding difference dimension polynomial is given
by$$\phi(r)=4r^4+\frac{56}{3}r^3+36r^2+4r+22.$$

Comparing difference dimension polynomials computed for the forward and symmetric difference schemes we can conclude that the strength of the system of difference equations obtained via forward difference scheme is higher than the strength of the system obtained with the use of symmetric difference scheme. This time we obtain  that the forward scheme is more preferable from the point of view of strength.
\end{example}

\begin{example}\textnormal{\textbf{(Equations for electromagnetic field given by potential)}}\\
 An electromagnetic field can be defined by the differential equations describing its potential, cf \cite[Ex. 9.2.6.]{KLMP}. The corresponding system of PDEs, which involves four unknown functions $\psi_1(x_1,\ldots,x_4),\ldots,\psi_4(x_1,\ldots,x_4)$, is as follows.
 \begin{eqnarray}\label{potential2}\sum_{j=1}^4\frac{\partial}{\partial x_j}\psi_j&=&0,\\\label{potential1}\sum_{j=1}^4\left(\frac{\partial^2}{\partial x_j^2}\psi_i-\frac{\partial^2}{\partial x_i~\partial x_j}\psi_j\right)&=&0.\end{eqnarray}

 \centerline{\textnormal{\bf{Differential dimension polynomial}}} Let $K$ be a differential field with basic set $\Delta=\{\delta_i=\D\frac{\partial}{\partial x_i}~|~i=1,\ldots,4\}$. Then equations \eqref{potential1} and \eqref{potential2} give rise to a differential $K[\Delta]$-module $M$ with generators $m_1,\ldots,m_4$ satisfying for $i=1,\ldots,4$ the defining equations\begin{eqnarray*}\sum_{j=1}^4\delta_j m_j&=&0,\\\sum_{j=1}^4\left(\delta_j^2 m_i-\delta_i\delta_j m_j\right)&=&0.\end{eqnarray*}
Then $M$ is isomorphic to the factor module of a free $K[\Delta]$-module with free generators $e_1,\ldots,e_4$ by its submodule $N$ generated by\begin{equation}\label{differentialpotential}\left\{\sum_{j=1}^4\delta_j e_j\right\}\cup\left\{\left.\sum_{j=1}^4\left(\delta_j^2 e_i-\delta_i\delta_j e_j\right)~\right|~i=1,\ldots,4\right\}.\end{equation}
 Defining an admissible order $\prec$ by\begin{eqnarray*}\lefteqn{\delta_1^{a_1}\delta_2^{a_2}\delta_3^{a_3}\delta_4^{a_4}e_{j_1}\prec\delta_1^{b_1}\delta_2^{b_2}\delta_3^{b_3}
 \delta_4^{b_4}e_{j_2}:\Longleftrightarrow}\nonumber\\&&(a_1+a_2+a_3+a_4,j_1,a_1,a_2,a_3,a_t)
 \nonumber\\&&<_{\lex}(b_1+b_2+b_3+b_4,j_2,b_1,b_2,b_3,b_4),\nonumber\end{eqnarray*}
we obtain the following Gr\"obner basis $G$ for $N$.
\begin{eqnarray*}G&=&\{\delta_{1}^2 e^{3}+\delta_{2}^2 e^{3}+\delta_{4}^2 e^{3}+\delta_{3}^2 e^{3}\nonumber\\&&~\delta_{1}^2 e^{2}+\delta_{2}^2 e^{2}+\delta_{4}^2 e^{2}+\delta_{3}^2 e^{2}\nonumber\\&&~\delta_{1}^2 e+\delta_{3}^2 e+\delta_{4}^2 e+\delta_{2}^2 e\nonumber\\&&~\delta_{1}e+\delta_{2} e^{2}+\delta_{3} e^{3}+\delta_{4} e^{4}\nonumber\\&&~\delta_{1}^2 e^{4}- \delta_{1} \delta_{4}e+\delta_{2}^2 e^{4}-\delta_{2} \delta_{4} e^{2}+\delta_{3}^2 e^{4}-\delta_{3} \delta_{4} e^{3}\}.\nonumber                                                                                                                                                                                                                                                                                                                                                                                                                                                                                                                                                                                                                                                                                                                                                                                                                                                                                                                                                                                                                                                                                                                                                                                                                                                                     \end{eqnarray*}
Applying Proposition \ref{prop1} we obtain that the differential dimension polynomial associated with \eqref{potential1} and \eqref{potential2} is given by$$\phi(r)=r^3+\frac{11}{2}r^2+\frac{17}{2}r+4.$$

\centerline{\textnormal{\bf{Difference dimension polynomial for forward difference scheme}}}
Let $K$ be an inversive difference field with basic set $\Delta=\{\delta_i:x_i\mapsto x_i+1~|~i=1,\ldots,4\}$. Replacing every occurrence of $\delta_k$ in (8) by $\delta_k-1$ ($k=1,\ldots,4$) we obtain that the desired dimension polynomial is the $\Delta^{\ast}$-dimension polynomial associated with the factor module of the free $K[\Delta^{\ast}]$-module $E= \D\sum_{i=1}^{4}K[\Delta^{\ast}]e_{i}$ by its $K[\Delta^{\ast}]$-submodule generated by the set \begin{eqnarray*}&&\left\{\sum_{j=1}^4\delta_j e_j-e_j\right\}\cup\nonumber\\&&\left\{\left.\sum_{j=1}^4\left(\delta_j^2 e_i-2\delta_je_i+e_i-\delta_i\delta_j e_j+\delta_ie_j+\delta_je_j-e_j\right)\right|i=1,\ldots,4\right\}.\end{eqnarray*}
Exploring the idea described at the end of subsection A, let us consider $K$ as a difference field with basic set$$\Sigma=\{\alpha_i:x_i\mapsto x_i+1,\beta_i:x_i\mapsto x_i-1~|~i=1,\ldots,4\}$$and let the admissible order $\prec$ be given by
\begin{eqnarray*}
 \lefteqn{\alpha_1^{a_1}\alpha_2^{a_2}\alpha_3^{a_3}\alpha_4^{a_4}\beta_1^{b_1}\beta_2^{b_2}\beta_3^{b_3}\beta_4^{b_4}e_i\prec\alpha_1^{c_1}\alpha_2^{c_2}\alpha_3^{c_3}\alpha_4^{c_4}\beta_1^{d_1}\beta_2^{d_2}\beta_3^{d_3}\beta_4^{d_4}e_j}\nonumber\\&:\Longleftrightarrow&(a_1+b_1+\cdots+a_4+b_4,i,a_1,\ldots,a_4,b_1,\ldots,b_4)\nonumber\\&&<_{\lex}(c_1+d_1+\cdots+c_4+d_4,i,c_1,\ldots,c_4,d_1,\ldots,d_4).
\end{eqnarray*}
Using the Maple package ``Ore\_Algebra'' \cite{OreAlgebra} for computing a Gr\"obner basis of the $K[\Sigma]$-submodule generated by
\begin{eqnarray*}&&\left\{\sum_{j=1}^4\alpha_j e_j-e_j\right\}\cup\nonumber\\&&\left\{\left.\sum_{j=1}^4\left(\alpha_j^2 e_i-2\alpha_je_i+e_i-\alpha_i\alpha_j e_j+\alpha_ie_j+\alpha_je_j-e_j\right)\right|i=1,\ldots,4\right\}\end{eqnarray*}we obtain the set of leading terms of the Gr\"obner basis
\begin{eqnarray*}
&\{\alpha_{4}\beta_{8}e_{1}, \alpha_{3}\beta_{7}e_{1}, \alpha_{2}\beta_{6}e_{1}, \alpha_{1}\beta_{5}e_{1}, \alpha_{1}^2e_{1}, \alpha_{4}\beta_{8}e_{2}, \alpha_{2}^2\beta_{5}e_{1},&\nonumber\\
& \alpha_{1}\beta_{6}\beta_{7}\beta_{8}^2e_{1},  \alpha_{3}\beta_{7}e_{2},\alpha_{2}\beta_{6}e_{2}, \alpha_{1}\beta_{5}e_{2}, \alpha_{1}^2e_{2}, \alpha_{4}\beta_{8}e_{3}, &\nonumber\\
&\beta_{8}e_{4}, \alpha_{4}\beta_{5}\beta_{6}\beta_{7}e_{1}, \alpha_{3}^2\beta_{5}\beta_{6}e_{1},  \alpha_{2}\beta_{5}^2\beta_{7}\beta_{8}^2e_{1},\alpha_{2}^2\beta_{5}e_{2},&\nonumber\\
&  \alpha_{1}\beta_{6}\beta_{7}\beta_{8}^2e_{2},\alpha_{3}\beta_{7}e_{3}, \alpha_{2}\beta_{6}e_{3}, \alpha_{1}\beta_{5}e_{3},  \alpha_{1}^2e_{3}, \alpha_{4}e_{4},\alpha_{3}\beta_{5}^2\beta_{6}^2\beta_{8}^2e_{1},&\nonumber\\
& \alpha_{4}\beta_{5}\beta_{6}\beta_{7}e_{2}, \alpha_{3}^2\beta_{5}\beta_{6}e_{2},  \alpha_{2}\beta_{5}^2\beta_{7}\beta_{8}^2e_{2},  \alpha_{2}^2\beta_{5}e_{3}, \alpha_{1}\beta_{6}\beta_{7}\beta_{8}^2e_{3}, \alpha_{3}\beta_{7}e_{4},&\nonumber\\
& \alpha_{2}\beta_{6}e_{4}, \alpha_{1}\beta_{5}e_{4}, \alpha_{1}^2e_{4}, \beta_{5}^2\beta_{6}^2\beta_{7}^2\beta_{8}^2e_{1}, \alpha_{3}\beta_{5}^2\beta_{6}^2\beta_{8}^2e_{2}, \alpha_{4}\beta_{5}\beta_{6}\beta_{7}e_{3},&\nonumber\\
&\alpha_{3}^2\beta_{5}\beta_{6}e_{3} , \alpha_{2}\beta_{5}^2\beta_{7}\beta_{8}^2e_{3}, \alpha_{2}^2\beta_{5}e_{4}, \alpha_{1}\beta_{6}\beta_{7}e_{4}, \beta_{5}^2\beta_{6}^2\beta_{7}^2\beta_{8}^2e_{2},&\nonumber\\
& \alpha_{3}\beta_{5}^2\beta_{6}^2\beta_{8}^2e_{3},  \alpha_{3}^2\beta_{5}\beta_{6}e_{4}, \alpha_{2}\beta_{5}^2\beta_{7}e_{4}, \beta_{5}^2\beta_{6}^2\beta_{7}^2\beta_{8}^2e_{3},&\nonumber\\
& \alpha_{3}\beta_{5}^2\beta_{6}^2e_{4}, \beta_{5}^2\beta_{6}^2\beta_{7}^2e_{4}\}.\nonumber
\end{eqnarray*}
Applying Proposition \ref{prop1} we compute the difference dimension polynomial associated with a forward difference scheme for \eqref{potential1} and \eqref{potential2} to be$$\phi(r)=15r^3-\frac{7}{2}r^2+\frac{43}{2}r^2+2.$$

\centerline{\textnormal{\bf{Difference dimension polynomial for symmetric difference scheme}}} ~\\Let $K$ be an inversive difference field with basic set $\Delta=\{\delta_i:x_i\mapsto x_i+1~|~i=1,\ldots,4\}$. Replacing every occurrence of $\delta_{k}$ in (8) by $\frac{1}{2}(\delta_k-\delta_k^{-1})$ ($k=1,\ldots,4$) we obtain the set
\begin{eqnarray*}&&\left\{\sum_{j=1}^4\frac{1}{2}(\delta_j -\delta_j^{-1})e_j\right\}\cup\nonumber\\&&\left\{\sum_{j=1}^4\frac{1}{4}\left(\delta_j^2 e_i-2e_i+\delta_j^{-2}e_i\right.\right.\nonumber\\
&&~\left.\left.\rule{0pt}{20pt} \left.-\delta_i\delta_j e_j+\delta_i\delta_j^{-1}e_j+\delta_i^{-1}\delta_je_j-\delta_i^{-1}\delta_j^{-1}e_j\right)~\right|~i=1,\ldots,4\right\},\end{eqnarray*}
which generates a $K[\Delta^{\ast}]$-submodule $N$ of the free $K[\Delta^{\ast}]$-module $M = \D\sum_{i=1}^{4}K[\Delta^{\ast}]e_{i}$ such that
the difference dimension polynomial of our system of difference equations is the $\Delta^{\ast}$-dimension polynomial of $M/N$.

Using the approach described at the end of subsection A, we treat $K$ as a difference field with basic set$$\Sigma=\{\alpha_i:x_i\mapsto x_i+1,\beta_i:x_i\mapsto x_i-1~|~i=1,\ldots,4\}$$and consider the admissible order $\prec$ given by
\begin{eqnarray*}
 \lefteqn{\alpha_1^{a_1}\alpha_2^{a_2}\alpha_3^{a_3}\alpha_4^{a_4}\beta_1^{b_1}\beta_2^{b_2}\beta_3^{b_3}\beta_4^{b_4}e_i\prec\alpha_1^{c_1}\alpha_2^{c_2}\alpha_3^{c_3}\alpha_4^{c_4}\beta_1^{d_1}\beta_2^{d_2}\beta_3^{d_3}\beta_4^{d_4}e_j}\nonumber\\&:\Longleftrightarrow&(a_1+b_1+\cdots+a_4+b_4,i,a_1,\ldots,a_4,b_1,\ldots,b_4)\nonumber\\&&<_{\lex}(c_1+d_1+\cdots+c_4+d_4,i,c_1,\ldots,c_4,d_1,\ldots,d_4).
\end{eqnarray*}
Once again using the Maple package ``Ore\_Algebra'' for computing a Gr\"obner basis of the $K[\Sigma]$-submodule generated by
\begin{eqnarray*}&&\left\{\sum_{j=1}^4\frac{1}{2}(\alpha_j -\beta_j)e_j\right\}\cup\nonumber\\&&\left\{\sum_{j=1}^4\frac{1}{4}\left(\alpha_j^2 e_i-2e_i+\alpha_j^{-2}e_i\right.\right.\nonumber\\
&&~\left.\left.\rule{0pt}{20pt} \left.\rule{0pt}{9pt}-\alpha_i\alpha_j e_j+\alpha_i\beta_je_j+\beta_i\alpha_je_j-\beta_i\beta_je_j\right)~\right|~i=1,\ldots,4\right\}\end{eqnarray*}we obtain the set of leading terms of the Gr\"obner basis
\begin{eqnarray*}
 &\{\alpha_{4}\beta_{8}e_{1}, \alpha_{3}\beta_{7}e_{1}, \alpha_{2}\beta_{6}e_{1}, \alpha_{1}\beta_{5}e_{1}, \alpha_{1}^2e_{1}, \alpha_{4}\beta_{8}e_{2}, \alpha_{2}^2\beta_{5}e_{1}, \alpha_{3}\beta_{7}e_{2},&\nonumber\\& \alpha_{2}\beta_{6}e_{2}, \alpha_{1}\beta_{5}e_{2}, \alpha_{1}^2e_{2}, \alpha_{4}\beta_{8}e_{3}, \beta_{8}^2e_{4}, \alpha_{3}^2\beta_{5}\beta_{6}e_{1}, \alpha_{2}^2\beta_{5}e_{2}, \alpha_{3}\beta_{7}e_{3},&\nonumber\\& \alpha_{2}\beta_{6}e_{3}, \alpha_{1}\beta_{5}e_{3}, \alpha_{1}^2e_{3}, \alpha_{4}e_{4}, \beta_{5}^3\beta_{6}\beta_{7}\beta_{8}e_{1}, \alpha_{4}^2\beta_{5}\beta_{6}\beta_{7}e_{1}, \alpha_{3}^2\beta_{5}\beta_{6}e_{2}, &\nonumber\\& \alpha_{2}^2\beta_{5}e_{3},\alpha_{3}\beta_{7}e_{4}, \alpha_{2}\beta_{6}e_{4}, \alpha_{1}\beta_{5}e_{4}, \alpha_{1}^2e_{4}, \beta_{5}^3\beta_{6}\beta_{7}\beta_{8}e_{2}, \alpha_{4}^2\beta_{5}\beta_{6}\beta_{7}e_{2},&\nonumber\\&   \alpha_{3}^2\beta_{5}\beta_{6}e_{3}, \alpha_{2}^2\beta_{5}e_{4},\beta_{5}^3\beta_{6}\beta_{7}\beta_{8}e_{3},\alpha_{4}^2\beta_{5}\beta_{6}\beta_{7}e_{3}, \alpha_{3}^2\beta_{5}\beta_{6}e_{4}, \beta_{5}^3\beta_{6}\beta_{7}e_{4}\}&\nonumber
\end{eqnarray*}
Applying Proposition \ref{prop1} we compute the difference dimension polynomial associated with a symmetric difference scheme for \eqref{potential1} and \eqref{potential2} to be$$\phi(r)=16r^3-8r^2+24r+8.$$
Comparing difference dimension polynomials computed for the forward and symmetric difference schemes we see that in this case, as in the previous example, the forward scheme is more preferable from the point of view of strength.
\end{example}

\section*{Conclusion}

We have developed a method for evaluation of the strength of systems of partial differential and difference equations based on the computation of the corresponding differential and difference dimension polynomials. We have also presented algorithms for such computation that extend the Gr\"obner basis technique to the cases of differential, difference, and inversive difference modules. Finally, we have determined the strength of some fundamental systems of PDEs of mathematical physics and the strength of the corresponding systems of partial difference equations obtained by the forward and symmetric difference schemes.

\section*{Acknowledgements}
The first author's research was partially supported by the Austrian Science Fund (FWF): W1214-N15, project DK11 and project no. P20336-N18 (DIFFOP) as well as the Austrian Marshall Plan Foundation: scholarship no. 256 420 24 7 2011.

\medskip

The second author's research was supported by the NSF Grant CCF 1016608

\end{document}